\documentclass[12pt,a4paper]{article}

\usepackage{amsmath,defns}
\usepackage[colour]{ajr}
\allowdisplaybreaks

\renewcommand{\includegraphics}[2][]{}

\newcommand{\KS}{Ku\-ra\-mo\-to--Siva\-shin\-sky}
\newcommand{\ibc}{\textsc{ibc}s}
\newcommand{\uj}{u_{j}}
\newcommand{\up}{u_{j+1}}
\newcommand{\um}{u_{j-1}}
\newcommand{\umm}{u_{j-2}}
\newcommand{\upp}{u_{j+2} }
\newcommand{\cM}{{\cal M}}
\newcommand{\cE}{{\cal E}}
\newcommand{\wt}[1]{\widetilde{#1}}
\newcommand{\DG}[3]{\frac{\partial^{#3} #2}{\partial #1^{#3}}}

\title{Accurately model the Kuramoto--Sivashinsky dynamics with
holistic discretisation}

\author{T. MacKenzie\thanks{Department of Mathematics and Computing,
University of Southern Queensland, Toowoomba, Queensland~4352,
\textsc{Australia}.} 
\and A.~J. Roberts\thanks{Department of Mathematics and Computing,
University of Southern Queensland, Toowoomba, Queensland~4352,
\textsc{Australia}.   \url{http://www.sci.usq.edu.au/staff/aroberts}}}

\date{March 23, 2005 --- without figures.\\ See
\url{http://www.sci.usq.edu.au/staff/aroberts/ksdoc.pdf} to download a
version with the figures}

\begin{document}

\maketitle

\begin{abstract}
We analyse the nonlinear \KS\ equation to develop accurate
discretisations modeling its dynamics on coarse grids.
The analysis is based upon centre manifold theory so we are assured
that the discretisation accurately models the dynamics and may be
constructed systematically.
The theory is applied after dividing the physical domain into small
elements by introducing isolating internal boundaries which are later
removed.
Comprehensive numerical solutions and simulations show that the
holistic discretisations excellently reproduce the steady states and the
dynamics of the \KS\ equation.
The \KS\ equation is used as an example to show how holistic
discretisation may be successfully applied to fourth order, nonlinear,
spatio-temporal dynamical systems.
This novel centre manifold approach is holistic in the sense that it
treats the dynamical equations as a whole, not just as the sum of
separate terms.
\end{abstract}


\paragraph{AMS Subj Class:} 37M99, 37L65, 65M20

\tableofcontents

\section{Introduction}

The \KS\ equation, here
\begin{equation} 
    \D tu+4\DDDD xu+\alpha\left(u\D xu+\DD xu\right)=0\,.
    \label{EpdeKS}
\end{equation}
was introduced by Sivashinsky~\cite{Sivas77} as a model of
instabilities on interfaces and flame fronts, and
Kuramoto~\cite{Kuramoto78} as a model of phase turbulence in chemical
oscillations.
It receives considerable attention as a model of complex
spatio-temporal dynamics~\cite[e.g.]{Hyman86b, Pomeau84, Cross93,
Holmes96}.
In the form~(\ref{EpdeKS}), with $2\pi$~periodic boundary conditions,
$\alpha$~is a bifurcation parameter that depends upon the size of the
typical pattern~\cite{Scovel88}.
The \KS\ equation includes the mechanisms of linear negative
diffusion~$\alpha u_{xx}$, high-order dissipation~$4u_{xxxx}$, and
nonlinear advection\slash steepening~$\alpha uu_x$.
The system~(\ref{EpdeKS}) has strong dissipative dynamics arising from
the fourth order dissipation.
Many modes of this system decay rapidly because of this strong
dissipation.
Thus the dynamics are dominated by a relatively few large scale modes.
We create and explore the macroscopic modelling of the \KS\ dynamics
using holistic discretisation as initiated by MacKenzie \&
Roberts~\cite{MacKenzie00a}.

We study the \KS\ equation here for several reasons.
Firstly, the \pde\ is fourth order and therefore, following the example
of Burgers' equation~\cite{Roberts98a}, provides a further test case
for the application of the holistic approach to higher order
dissipative \pde s.
Secondly, the \KS\ equation is analogous to the Navier--Stokes
equations of fluid dynamics.
Holmes, Lumley \& Berkooz~\cite{Holmes96} argued that these analogies
exist on two levels: in the energy source and dissipation terms of both
dynamical systems; and in the reflection and translational symmetries
of the \KS\ equation and the spanwise symmetries of the Navier--Stokes
equations in the boundary layer.
This analogy between symmetries suggests the Fourier series and
corresponding modal interactions are comparable for these two problems.
Thirdly, Cross \& Hohenberg~\cite{Cross93} describes how the \KS\
equation exhibits the complexities of weak turbulence or
spatio-temporal chaos.
The complex dynamics of the \KS\ equation~(\ref{EpdeKS}) is searching
test of the performance of the holistic approach to coarse grained
modelling of dynamical systems.

Approximate inertial manifolds and variants \cite[e.g.]{Foias85b,
Foias88c, Foias94, Armbruster89, Jolly90} capture the long-term low
dimensional behaviour of the \KS\ equation.
Most constructions of approximate inertial manifolds are based upon
nonlinear Galerkin methods~\cite[e.g.]{Roberts89, marion89, Jolly90,
Foias94}.
Approximate inertial manifolds are generally constructed by finding
\emph{global} eigenfunctions of the linear dynamics.
Our approach is similar to these methods in that we project onto
natural solutions of the \pde, and performs nearly as well, see
\S\ref{S_ks_ss_g}.
But in contrast, the holistic approach undertaken here bases analysis
upon the \emph{local} dynamics within and between finite elements and
thus we contend it will be more useful in applications; for example, it
is readily adapted to the modelling of a wide variety of physical
boundary conditions~\cite{Roberts01b}.

Our approach is to divide the spatial domain into disjoint finite
elements (\S\ref{SSibc}).
Initially these finite elements are decoupled and so dissipation causes
the solution to exponentially quickly become constant in each element.
We then couple the elements together so that information is exchanged
between elements---parameterised by a coupling parameter~$\gamma$ so
that $\gamma=1$ recovers the original \KS\ equation.
The coupling drives the evolution of the field in each element.
We solve the \KS\ \pde\ within each element with the coupling and hence
resolve subgrid scale dynamics.
Centre manifold theory \cite[e.g.]{Carr81, Roberts97a} then provides
the rigorous support for holistic models as introduced by
Roberts~\cite{Roberts98a} for Burgers' equation and discussed
in~\S\ref{S_KScm}.

A low order analysis, reported in \S\ref{SSsome}, of the \KS\
equation~(\ref{EpdeKS}) favours the discretisation
\begin{eqnarray}
&&
    \frac{d\uj}{dt}+\frac{4\upp-16\up+24\uj-16\um+4\umm}{h^4}
    \nonumber \\ &&
    {}
    +\alpha\left(
    \frac{-\upp+16\up-30\uj+16\um-\umm}{12h^2}\right)
    \nonumber \\ &&
    {}+\alpha\left(u_j\frac{\up-\um}{4h}+\frac{\up^2-\um^2}{4h}
    -\frac{\upp\up-\umm\um}{12h}
    \right)
    \approx0
    \,,
    \label{Eks_fd}
\end{eqnarray}
where the $u_j$s~are grid values spaced $h$~apart.
The first two lines of the holistic discretisation~(\ref{Eks_fd}) shows the
holistic method generates conventional centered finite difference
approximations for the linear terms $4u_{xxxx}$~and~$\alpha u_{xx}$.
The third line details a specific nonstandard approximation for the
nonlinear term~$\alpha uu_x$: it is a mix of three valid approximations
to~$uu_x$; the specific mix is determined by the subgrid scale
modelling of physical processes in the holistic
approach,  see~\S\ref{S_ks_sg}.
The holistic discretisation is not constructed by discretising the \KS\
equation~(\ref{EpdeKS}) term by term, rather the subgrid scale dynamics
of~(\ref{EpdeKS}) together with inter-element coupling generate the
specific holistic discretisation~(\ref{Eks_fd}).

The discretisation~(\ref{Eks_fd}) is a low-order approximation.
Centre manifold theory provides systematic refinements.
Analysis to higher orders in nonlinearity or inter-element interaction,
discussed in~\S\ref{S_ks_rel}, gives further refinement to the
discretisation.
The higher order terms come from resolving more subgrid scale
interactions.
These higher order analyses lead to higher order consistency, as
element size~$h\to0$\,, between the equivalent \pde s of the holistic
discretisations, such as (\ref{Eks_fd}), and the \KS{} \pde\
(see~\S\ref{S_ks_epde}).
Such consistency is further justification for our approach.

The bulk of this paper is then a comprehensive comparative study of the
various models of the \KS\ dynamics; further details are reported
by MacKenzie~\cite{MacKenzie05}.
A detailed numerical study of the holistic predictions for the steady
states of the \KS\ equation is the focal point of
Section~\ref{chap_ks_ss}, followed by an exploration of the holistic
predictions for the time dependent phenomena of the \KS\ equation in
Section~\ref{chap_ks_td}.
We look at: the predicted steady states, their stability and compare
bifurcation diagrams; the dynamics near the steady states; Hopf
bifurcations leading to period doubling sequences; and the
spatio-temporal patterns at relatively large nonlinearity
parameter~$\alpha$.
We find that the holistic models have excellent performance on coarse
grids thus leading to simulations that may use large time steps.
The excellent performance detailed herein is further
evidence that the holistic approach is a robust and useful method for
discretising \pde{}s.

\section{Use a homotopy in the inter-element coupling}
\label{S_KSho}

The construction of a discretisation is based upon breaking the spatial
domain into disjoint finite elements and then joining them together
again.
We control this process by a coupling parameter~$\gamma$ that
smoothly parametrises the transition between decoupled elements and
fully coupled elements for which we recover a model
for the original \pde.
Furthermore, we construct the model using solutions of the \pde\ within
each element and hence resolve subgrid scale dynamics.
Centre manifold theory \cite[e.g.]{Carr81, Roberts97a} provides the
rigorous support for holistic models as introduced by
Roberts~\cite{Roberts98a} for Burgers' equation.

\subsection{Introduce internal boundaries between elements}
\label{SSibc}

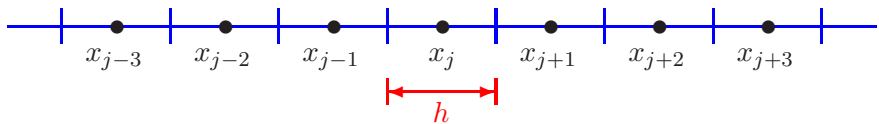
\begin{figure}
\begin{center}
\small \setlength{\unitlength}{0.3em}
\begin{picture}(80,15)
\thicklines \put(-10,10){\color{blue} \line(1,0){100}}
\multiput(-3.75,8)(12.5,0){8}{\color{blue} \line(0,1){4}}
\multiput(2.65,10)(12.5,0){7}{\circle*{1.5}} \put(-1,6){$x_{j-3}$}
\put(11.5,6){$x_{j-2}$} \put(24,6){$x_{j-1}$} \put(38.5,6){$x_j$}
\put(49,6){$x_{j+1}$} \put(61.5,6){$x_{j+2}$}
\put(74,6){$x_{j+3}$} \multiput(33.75,1)(12.5,0){2}{\color{red}
\line(0,1){3}} \put(33.75,2.5){\color{red} \vector(1,0){12.5}}
\put(46.25,2.5){\color{red} \vector(-1,0){12.5}}
\put(39,-1){\color{red} $h$}
\end{picture}
\end{center}
\caption{An example of the $1D$ grid with regular elements of
width~$h$. The $j$th~element is centered about the grid point~$x_j$.
The vertical blue lines form the element boundaries, which for the
$j$th element are located at $x_{j\pm1/2}=(j\pm1/2)h$.}
\label{fig_ks_grid}
\end{figure}

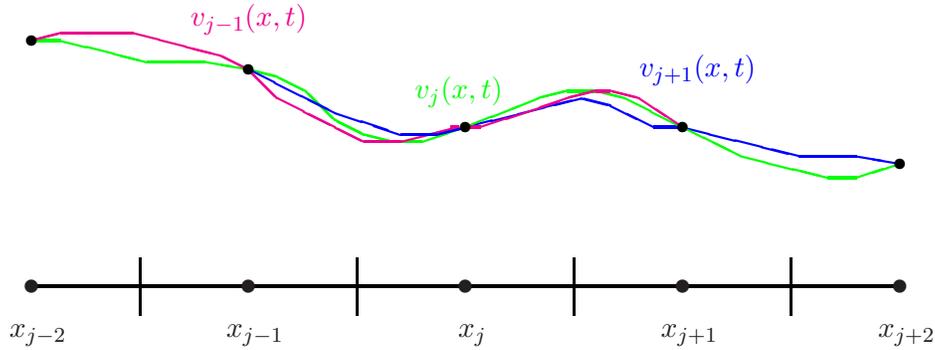
\begin{figure}
\begin{center}
\small \setlength{\unitlength}{0.25em}
\begin{picture}(80,50)
\thicklines
\put(-20,10){\color{black} \line(1,0){120}}
\multiput(-5,6)(30,0){4}{\color{black} \line(0,1){8}}
\multiput(-20,10)(30,0){5}{\circle*{2}} \put(7,3){{$x_{j-1}$}}
\put(39,3){{$x_j$}} \put(67,3){{$x_{j+1}$}}
\put(-23,3){{$x_{j-2}$}} \put(97,3){{$x_{j+2}$}}
\put(2,46){{\color{magenta} $v_{j-1}(x,t)$}}
\put(-20,44){\color{green} \line(1,0){4}}
\put(-16,44){\color{green} \line(4,-1){12}}
\put(-4,41){\color{green} \line(1,0){8}}
\put(4,41){\color{green} \line(6,-1){6}}
\put(10,40){\color{green} \line(4,-1){4}}
\put(14,39){\color{green} \line(2,-1){4}}
\put(18,37){\color{green} \line(1,-1){4}}
\put(22,33){\color{green} \line(2,-1){4}}
\put(26,31){\color{green} \line(4,-1){4}}
\put(30,30){\color{green} \line(1,0){4}}
\put(34,30){\color{green} \line(3,1){6}}
\put(-20,44){\color{magenta} \line(4,1){4}}
\put(-16,45){\color{magenta} \line(1,0){10}}
\put(-6,45){\color{magenta} \line(4,-1){12}}
\put(6,42){\color{magenta} \line(2,-1){4}}
\put(10,40){\color{magenta} \line(1,-1){4}}
\put(14,36){\color{magenta} \line(2,-1){8}}
\put(22,32){\color{magenta} \line(2,-1){4}}
\put(26,30){\color{magenta} \line(1,0){5}}
\put(31,30){\color{magenta} \line(4,1){8}}
\put(38,32){\color{magenta} \line(1,0){1}}
%\put(-20,40){\color{blue} \line(6,-1){6}}
%\put(-14,39){\color{blue} \line(5,1){10}}
%\put(-4,41){\color{blue} \line(6,1){6}}
%\put(2,42){\color{blue} \line(1,0){4}}
%\put(6,42){\color{blue} \line(2,-1){4}}
\put(10,40){\color{blue} \line(2,-1){12}}
\put(22,34){\color{blue} \line(3,-1){9}}
\put(31,31){\color{blue} \line(1,0){5}}
\put(36,31){\color{blue} \line(4,1){4}}
\put(40,32){\color{blue} \line(4,1){16}}
\put(56,36){\color{blue} \line(4,-1){4}}
\put(60,35){\color{blue} \line(2,-1){6}}
\put(66,32){\color{blue} \line(1,0){4}}
\put(70,32){\color{blue} \line(4,-1){16}}
\put(86,28){\color{blue} \line(1,0){8}}
\put(94,28){\color{blue} \line(6,-1){6}}
\put(33,36){{\color{green} $v_{j}(x,t)$}}
\put(40,32){\color{green} \line(5,2){10}}
\put(50,36){\color{green} \line(4,1){4}}
\put(54,37){\color{green} \line(1,0){4}}
\put(58,37){\color{green} \line(4,-1){4}}
\put(62,36){\color{green} \line(2,-1){16}}
\put(78,28){\color{green} \line(4,-1){12}}
\put(90,25){\color{green} \line(1,0){4}}
\put(94,25){\color{green} \line(3,1){6}}
\put(40,32){\color{magenta} \line(1,0){2}}
\put(42,32){\color{magenta} \line(3,1){12}}
\put(54,36){\color{magenta} \line(4,1){4}}
\put(58,37){\color{magenta} \line(1,0){2}}
\put(60,37){\color{magenta} \line(4,-1){4}}
\put(64,36){\color{magenta} \line(3,-2){6}}
%\put(70,32){\color{magenta} \line(3,-2){6}}
%\put(76,28){\color{magenta} \line(2,-1){8}}
%\put(84,24){\color{magenta} \line(4,-1){8}}
%\put(92,22){\color{magenta} \line(1,0){4}}
%\put(96,22){\color{magenta} \line(4,1){4}}
\put(-20,44){\color{black} \circle*{1.5}}
\put(10,40){\color{black} \circle*{1.5}}
\put(40,32){\color{black} \circle*{1.5}}
\put(70,32){\color{black} \circle*{1.5}}
\put(100,27){\color{black} \circle*{1.5}}
\put(64,39){{\color{blue} $v_{j+1}(x,t)$}}
\end{picture}
\end{center}
\caption{Schematic diagram of the fields {\color{green} $v_j(x,t)$},
{\color{blue} $v_{j+1}(x,t)$} and {\color{magenta} $v_{j-1}(x,t)$} for
the non-local \ibc~(\ref{EbcsddKS}--\ref{EsbcidKS}) with $\gamma=1$\,.
See the fields pass through neighbouring grid values $u_j$
and~$u_{j\pm1}$, and also $u_{j\pm2}$ when appropriate.}
\label{fig_nl_ibc}
\end{figure}

Establish the spatial discretisation by dividing the domain into
$m$~elements of equal and finite width~$h$ and introducing an
equispaced grid of collocation points, $x_j=jh$, at the centre of each
element, see Figure~\ref{fig_ks_grid}.\footnote{In principle, elements
may be of unequal size.
However, to simplify the analysis, herein all elements will be of equal
width~$h$.} Express the subgrid field in the $j$th~element by
$u=v_j(x,t)$ --- we solve the \KS\ \pde~(\ref{EpdeKS}) with
inter-element coupling introduced via artificial internal boundary
conditions~(\ibc).
We introduce a homotopy in an inter-element coupling parameter
$\gamma$: when $\gamma=0$ the elements are effectively isolated from
each other, providing the basis for the application of centre manifold
theory; whereas when evaluated at $\gamma=1$, the elements are fully
coupled together and hence the discretised model applies to the
original \pde.
Since the \KS\ equation is fourth order we require four \ibc\ for each
element to ensure satisfactory coupling between neighbouring elements.
Here we use the non-local \ibc{}
\begin{eqnarray}
 \delta_x v_j(x,t)&=&\gamma\,\delta v_{j\pm1/2}(x,t)
\quad\mbox{at }{x=x_{j\pm1/2}}\,,\label{EbcsddKS}\\
\delta^3_x v_j(x,t)&=&\gamma^2\delta^3 v_{j\pm1/2}(x,t)
\quad\mbox{at }{x=x_{j\pm1/2}}\,,\label{EsbcidKS}
\end{eqnarray}
which are an extension of the non-local \ibc\ explored by
Roberts~\cite{Roberts00a} for Burgers' equation; a local possibility 
for the \ibc{}s was explored by MacKenzie~\cite{MacKenzie05}.
These non-local \ibc\ involve the centered difference operators
$\delta$~and~$\delta_x$: the operator~$\delta_x$ denotes a centered
difference in~$x$ only, with step~$h$; whereas the operator~$\delta$
denotes a centered difference applied to the grid index~$j$ with
step~1; so for example, the first \ibc~(\ref{EbcsddKS}) is
\begin{equation}
 v_j(x_{j\pm1},t)-v_j(x_j,t)
 =\gamma[ v_{j+1}(x_{j\pm1},t)-v_j(x_j,t) ]
\,.
\end{equation}
Note: the field $v_j(x,t)$ extends analytically to at least
$x_{j\pm2}$, to allow the application of the non-local
\ibc~(\ref{EsbcidKS}).
The physical interpretation of these \ibc\ is not obvious.
Firstly, when $\gamma=0$, (\ref{EbcsddKS}--\ref{EsbcidKS}) ensures the
first and third differences in~$x$ of the field~$v_j$ centered about
the element boundaries $x_{j\pm1/2}$ are zero.
These isolate each element from its neighbours as there is then no
coupling between them.
In each element $v_j(x,t)=\mbox{constant}$ is an equilibrium.
It is dynamically attractive provided the instability controlled by
$\alpha/h^2$ is not too large compared with the dissipation of
order~$1/h^4$.
This simple class of piecewise constant solutions provide the basis for
analysing the $\gamma\ne0$ case when the elements are coupled together.
Secondly, the non-local \ibc\ evaluated at $\gamma=1$ requires that the
field~$v_{j}(x,t)$, when extrapolated to $x_{j\pm1}$~and~$x_{j\pm2}$,
is to equal the grid point value of the subgrid field of that element,
$u_{j\pm1}$~and~$u_{j\pm2}$ respectively.
See the schematic representation in Figure~\ref{fig_nl_ibc} of these
non-local boundary conditions evaluated at $\gamma=1$.
This restores sufficient continuity to ensure the holistic model
applies to the original \pde.

The inter-element coupling parameter~$\gamma$ controls the flow of
information between neighbouring elements.
We construct solutions as power series expansions in the coupling
parameter~$\gamma$.\footnote{Such homotopies are used successfully in
other numerical methods.
For example, Liao~\cite{Liao95} proposed a homotopy in his general
boundary element method from auxiliary linear operators whose
fundamental solutions are well known.
In our application the homotopy is only in the \ibc.} When
$\Ord{\gamma^2}$ terms are neglected in the holistic model, the field
in the $j$th~element involves information about the fields in the
$j\pm1$ elements.
Similarly, when $\Ord{\gamma^3}$ terms are neglected in the
approximation, the field in the $j$th~element involves information
about the fields in the $j\pm1$~and~$j\pm2$ elements.
Consequently, the order of~$\gamma$ retained in the holistic
model controls the stencil width of the discretisation.

Roberts~\cite{Roberts00a} argued that this particular form of the
non-local \ibc\ ensure that these holistic models are consistent with
any given \pde\ to high orders in the grid size~$h$ as $h\to 0$\,.

\subsection{Centre manifold theory supports the discretisation}
\label{S_KScm}

The existence, relevance and approximation theorems~\cite[e.g.]{Carr81,
Carr83b} of centre manifold theory apply to the \KS\
\pde~(\ref{EpdeKS}) with \ibc~(\ref{EbcsddKS}--\ref{EsbcidKS}).
Similar to the application to Burgers' equation by
Roberts~\cite{Roberts98a}, the result here is support for a low
dimensional discrete model for the \KS\ dynamics at finite grid size.

% \subsubsection{There exists a centre manifold}
% \label{S_cm_ex_ks}

Theoretical support is based upon the piecewise constant solutions
obtained when all the elements are insulated from each other.
Adjoin to the \KS\ \pde~(\ref{EpdeKS}) the dynamically trivial
equations for the coupling parameter~$\gamma$ and the nonlinearity
parameter~$\alpha$,
\begin{equation}
    \D t \gamma=\D t \alpha=0\,,
    \label{EtrivKS}
\end{equation}
and consider the dynamics in the extended state space
$(u(x),\gamma,\alpha)$.
Adjoining such trivial equations for parameters is commonly used to
unfold bifurcations \cite[\S1.5]{Carr81}.
In this extended space there is a subspace of fixed points with
$u=\mbox{constant}$ in each element and $\gamma=\alpha=0$\,.
Linearizing the \pde\ and \ibc\ about each fixed point,
$u=\mbox{constant}+u'(x,t)$\,, gives
\begin{displaymath}\quad
    \frac{\partial u'}{\partial t}=-\frac{\partial^4 u'}{\partial x^4}\,,
    \quad\mbox{such that}\quad
    \left.\delta_x u'(x,t)\right|_{x=x_{j\pm 1/2}}=
    \left.\delta^3_x u'(x,t)\right|_{x=x_{j\pm 1/2}}=0\,.
\end{displaymath}
The $n$th~linear eigenmode associated with each element is
\begin{equation}
\alpha=\gamma=0\,,\quad
    u'\propto
        e^{\lambda_n t}\cos\left[\frac{n\pi}{h}(x-x_{j-1/2})\right]\,,
    \label{EmodeKSl}
\end{equation}
for the non-local \ibc~(\ref{EbcsddKS}--\ref{EsbcidKS}), where
$n=0,1,\ldots$ and the eigenvalue $\lambda_n=-n^4\pi^4/h^4$\,.
There are also the trivial modes $\gamma=\mbox{const}$ and
$\alpha=\mbox{const}$.
Therefore, in a spatial domain of $m$~elements there are $m+2$ zero
eigenvalues: one associated with each of the $m$~elements; and two from
the trivial~(\ref{EtrivKS}).
All other eigenvalues are negative and $\le-\pi^4/h^4$.
Thus, the existence theorem, see~\cite[p.281]{Carr83b}
or~\cite[p.96]{Vanderbauwhede89}, guarantees that a $m+2$~dimensional
centre manifold~$\cM$ exists for the \KS\ \pde~(\ref{EpdeKS}) with the
trivial~(\ref{EtrivKS}) and \ibc~(\ref{EbcsddKS}--\ref{EsbcidKS}).

We parametrise the $m+2$ dimensional centre manifold~$\cM$ by the $m+2$
parameters $\gamma$, $\alpha$ and the grid values~$\uj$.\footnote{These
grid values are one choice for the measure of the field~$u$ in each
element.
Other choices are possible, but the grid values appear most
convenient.} Denote $\vec u$ as the vector of the $m$~grid values.
Thus for some function~$v$ the centre manifold~$\cM$ is
\begin{equation}
    u(x,t)=v(x;\vec u,\gamma,\alpha)\,.
    \label{EcmvKS}
\end{equation}
Equivalently, we decompose the centre manifold field into that for
each element: 
\begin{equation}
u=v(x;\vec{u},\gamma,\alpha)
=\sum_j v_j(x;\vec{u},\gamma,\alpha)\,\chi_j(x)\,,
\label{E_ks_cm_ab}
\end{equation}
where the characteristic function $\chi_j(x)$ is~1 when
$x_{j-1/2}<x<x_{j+1/2}$\,, and 0~otherwise; view the centre manifold as
the union of all states of the collection of subgrid fields
$v_j(x;\vec{u},\gamma,\alpha)$ over the physical domain.
The corresponding amplitude condition, that the field in each element
has to pass through its grid value, is
\begin{equation}
\uj=v(x_j;\vec u,\gamma,\alpha)\,.\label{E_ks_amp}
\end{equation}
The existence theorem~\cite{Carr83b} also asserts that on the centre
manifold the grid values~$\uj$ evolve deterministically in time
according to the system of \ode{}s
\begin{equation}
    \dot\uj=d\uj/dt=g_j(\vec u,\gamma,\alpha)\,,
    \label{EcmgKS}
\end{equation}
where $g_j$~is the restriction of the \KS\
\pde~(\ref{EpdeKS}) with the trivial~(\ref{EtrivKS}) and
\ibc~(\ref{EbcsddKS}--\ref{EsbcidKS}) to the centre manifold~$\cM$.
It is this evolution~(\ref{EcmgKS}) of the grid values that gives the
holistic discretisation.

Note that the centre manifold~$\cM$ is global in~$u$ but local in
$\gamma$~and~$\alpha$.
When the parameters $\gamma=\alpha=0$ the \KS\ \pde\ has a
$m$~dimensional ``centre'' subspace~$\cE$ of fixed points with the
field~$u$ being independently constant in each element; these are fixed
points for all~$\vec u$.
When the parameters $\gamma$~and~$\alpha$ are non-zero this subspace is
``bent'' to the curved centre manifold~$\cM$.
Thus the models we construct are valid for small enough
$\gamma$~and~$\alpha$, although we use them at finite
$\gamma$~and~$\alpha$, but are formally valid for all~$|\vec u|$.
Numerical solutions of the centre manifold models, such as those in
\S\ref{S_ks_sg}, indicate that parameter values as large as $\gamma=1$
and $\alpha=20$--$50$ are indeed within the range of validity of our
approach, even on relatively coarse grids.

We now support the claim that the evolution of the discrete grid
values~(\ref{EcmgKS}) actually models the \KS\ system~(\ref{EpdeKS}).
The relevance theorem of centre manifolds, \cite[p.282]{Carr83b}
or~\cite[p.128]{Vanderbauwhede89}, guarantees that all solutions of the
\KS\ system~(\ref{EpdeKS}) with (\ref{EtrivKS}) and the
\ibc~(\ref{EbcsddKS}--\ref{EsbcidKS}), which remain in some
neighbourhood of the subspace~$\cE$ in $(u(x),\gamma,\alpha)$ space are
exponentially quickly attracted to the centre manifold~$\cM$ and thence
to a solution of the $m$~discrete \ode{}s~(\ref{EcmgKS}).
For our application of centre manifold theory to the holistic model we
seek regimes where this neighbourhood includes $\gamma=1$~and~$\alpha$
of interest.
We estimate the rate of attraction by the leading negative eigenvalue,
here $\lambda_1=-\pi^4/h^4$\,.
The actual rate of attraction may be less due to the difference between
centre manifold~$\cM$ and centre subspace~$\cE$, but $\lambda_1$ will
be the correct order of magnitude.
This ensures the so-called asymptotic completeness~\cite{Robinson96}:
after the exponentially quick transients of the approach to~$\cM$ by
any trajectory, the evolution of the discretisation~(\ref{EcmgKS})
on~$\cM$ accurately models the \KS\ \pde~(\ref{EpdeKS}).

\subsection{Approximate the shape of the centre manifold}
\label{S_ks_sym}

Having established that we may find a low dimensional
description~(\ref{EcmvKS}--\ref{EcmgKS}) of the interacting elements
that is relevant to the \KS\ system~(\ref{EpdeKS}), we need to
construct the shape of centre manifold and the corresponding evolution
on the manifold.

The approximation theorem of Carr \& Muncaster \cite[p.283]{Carr83b}
assures us that upon substituting the
ansatz~(\ref{EcmvKS}--\ref{EcmgKS}) into the complete system and
solving to some order of error in $\alpha$~and~$\gamma$, then $\cM$ and
the evolution thereon will be approximated to the same order.
However, we need to evaluate the approximations at the coupling
parameter $\gamma=1$ because it is only then that the artificial
internal boundaries are removed.
Thus the actual error of the model due to the evaluation at $\gamma=1$
is not estimated.
However, the holistic method for discretising the \KS\ equation is
supported three ways: firstly, the smooth homotopy from $\gamma=0$ with
large spectral gap to the gravest decaying mode with decay rate
$\approx -\pi^4/h^4$; secondly the holistic models are consistent with
the \KS\ \pde to high order in grid size~$h$, see~\S\ref{S_ks_epde};
thirdly, we see in Sections~\ref{chap_ks_ss}--\ref{chap_ks_td} the
holistic models model accurately both steady state solutions and time
dependent phenomena of the \KS\ system.

To construct the centre manifold, we solve for the field~$v_j$ in each
element.
For definiteness, here we consider domains periodic in space, or
equivalently elements far from the influence of any physical boundary.
By translational symmetry of the \KS\ \pde~(\ref{EpdeKS}) the subgrid
field in each element is identical, except for the appropriate shift in
the grid index~$j$.
Thus, we construct the subgrid field and evolution for a general
$j$th~element, see some examples in Section~{\ref{S_ks_rel}}

The algebraic details of the derivation of the centre manifold
model~(\ref{EcmvKS}--\ref{EcmgKS}) are handled by computer algebra.
In an algorithm introduced by Roberts~\cite{Roberts96a}, iteration
drives to zero the residuals of the governing \pde~(\ref{EpdeKS}) and
its \ibc~(\ref{EbcsddKS}--\ref{EsbcidKS}) and amplitude
condition~(\ref{E_ks_amp}).
Since the algebraic details of the construction are tedious, they are
not given;  instead see the computer algebra procedure of~\cite{Roberts02a}.

This computer algebra is based upon driving the residuals of the
governing equations to zero in the following manner.
Recall from \S\ref{S_KScm} that the centre manifold~(\ref{EcmvKS}) is
parametrised by the grid values~$\vec{u}$ and that the evolution of the
grid values is given by~(\ref{EcmgKS}).
Thus substitute these into the \KS\ \pde~(\ref{EpdeKS}) and seek to
solve
\begin{equation}
\D {t}{v_j}=
\sum_k\frac{\partial
v_j}{\partial
u_{k}} g_{k}=-4\DDDD {x}{v_j}-\alpha\left(\DD {x}{v_j}
+v_j\D {x}{v_j}\right)\,,
\label{E_ks_symeq}
\end{equation}
together with the non-local \ibc~(\ref{EbcsddKS}--\ref{EsbcidKS}) and
the amplitude equation~(\ref{E_ks_amp}), to some order in parameters
$\gamma$~and~$\alpha$.
The iteration is that given any approximation, denoted by $\wt\ $, we
seek corrections, denoted by primes, such that $v_j=\wt{v}_j+v_j'$ and
$g_j=\wt{g}_j+g_j'$, better satisfy the \KS\ \pde.
Thus in each iteration we solve a problem of the form,
\begin{equation}
-4\DDDD x{v'_j}=g_j'+\mbox{Residual}\,,
\label{E_ks_iteq}
\end{equation}
where
\begin{equation}
\mbox{Residual}=\sum_k\frac{\partial
\wt{v}_j}{\partial
u_{k}} g_{k}+4\DDDD {x}{\wt{v}_j}
+\alpha\left(\DD {x}{\wt{v}_j}+\wt{v}_j\D {x}{\wt{v}_j}\right)\,,
\label{E_ks_reseq}
\end{equation}
together with the \ibc{}, for the corrections, primed quantities, to
the subgrid field and the evolution of the grid values.
Note: the residual in~(\ref{E_ks_reseq}) is the residual of the \KS\
system for the current approximation.
The iteration scheme starts with the linear solution in each element,
namely $v_j(x,\vec{u},\gamma,\alpha)=u_j$ and
$g_j(\vec{u},\gamma,\alpha)=0\,$.
The iteration terminates when the residuals of the \KS\
\pde~(\ref{E_ks_symeq}), and the \ibc, are zero to some order in
$\gamma$~and~$\alpha$.
Then theory assures us that the subgrid field in each element and the
evolution of the grid values are correct to the same order in
$\gamma$~and~$\alpha$.

\section{Various holistic models}
\label{S_ks_rel}

Here we record holistic models of the \KS\ \pde~(\ref{EpdeKS}), to
various orders in coupling parameter~$\gamma$, governing the width of
the numerical stencil, and in the nonlinearity parameter~$\alpha$.
For use the models need to be evaluated at $\gamma=1$ as then the
non-local \ibc~(\ref{EbcsddKS}--\ref{EsbcidKS}) ensure sufficient
continuity in the solution field.
We write the models in terms of the centered difference and mean
operators,
\begin{displaymath}
    \delta u_j=u_{j+1/2}-u_{j-1/2}
    \quad\text{and}\quad
    \mu u_j=(u_{j+1/2}+u_{j-1/2})/2\,,
\end{displaymath}
respectively.
The models are constructed using a \textsc{reduce} program
adapted from~\cite{Roberts02a}.
We only present in detail here holistic models to errors
$\Ord{\alpha^2}$ as the level of complexity increases enormously with
the order of~$\alpha$.

\subsection{Some holistic discretisations}
\label{SSsome}

In order to represent the spatial fourth derivative in the \KS{}
equation, we need at least a 5~point stencil approximation.
Thus we determine the interactions between at least
next-nearest neighbouring elements by obtaining up to at least
quadratic terms in the coupling parameter~$\gamma$.

\paragraph{The $\Ord{\gamma^3,\alpha^2}$ holistic discretisation} is
\begin{eqnarray}
\dot{u_j}&=&-\frac{\gamma \alpha}{h^2}\delta^2 u_j
-\frac{\gamma \alpha}{h}u_j\delta\mu u_j
-\frac{4\gamma^2}{h^4}\delta^4 u_j
+\frac{\gamma^2 \alpha}{12h^2}\delta^4 u_j
\nonumber\\&&{}+\frac{\gamma^2 \alpha}{12h}\left(
2u_j\delta^3\mu u_j+\delta^2 u_j
\delta^3\mu u_j +\delta^4 u_j \delta \mu u_j
\right)
\nonumber\\&&{}+{\cal O}\left(\gamma^3,\alpha^2\right)\,,
\label{E_hol8dg1r1}
\end{eqnarray}
for the non-local \ibc~(\ref{EbcsddKS}--\ref{EsbcidKS}).
This forms a basic 5~point stencil approximation, since the
evolution~$\dot u_j$ involves just~$u_j$, $u_{j\pm1}$~and~$u_{j\pm2}$.
The first line of~(\ref{E_hol8dg1r1}), when evaluated at $\gamma=1$,
gives a 2nd~order centered difference approximation for the
hyperdiffusion term~$4u_{xxxx}$, a 4th~order centered difference
approximation to the linear growth term~$\alpha u_{xx}$, and a
2nd~order centered difference approximation to the nonlinear advection
term~$\alpha uu_x$.
The second line modifies the nonlinear discretisation to account for
interaction with effects caused by the next-nearest neighbour elements.

The holistic discretisation~(\ref{E_hol8dg1r1})
contains the approximation
\begin{equation}
 uu_x|_{x_j} \approx \left(u_j\frac{\up-\um}{4h}
+\frac{\up^2-\um^2}{4h}-\frac{\upp\up-\umm\um}{12h}\right)\,.
\label{E_ks_hol_non2}
\end{equation}
when evaluated at $\gamma=1$.
This is a $1/2:1:-1/2$ mix of the approximations
\begin{equation}
\left. uu_x\right|_{x_j}
\approx u_j\frac{\up-\um}{2h}
\approx \frac{\up^2-\um^2}{4h}
\approx \frac{u_{j+2}u_{j+1}-u_{j-2}u_{j-1}}{6h}\,,
\label{Euux4}
\end{equation}
respectively.  This particular nonstandard
approximation~(\ref{E_ks_hol_non2}) to the nonlinear term~$\alpha
uu_x$, arises due to the modelling of subgrid scale interactions
between the \KS\ equation and the inter-element coupling.  Such
nonstandard approximations generated through this approach can have
robust numerical characteristics~\cite[\S2]{Roberts00a}.

\paragraph{The $\Ord{\gamma^4,\alpha^2}$ holistic discretisation} is
\begin{eqnarray}
\dot{u_j}&=&-\frac{\gamma \alpha}{h^2}\delta^2 u_j -\frac{\gamma
\alpha}{h}u_j\delta\mu u_j -\frac{4\gamma^2}{h^4}\delta^4 u_j
+\frac{\gamma^2 \alpha}{12h^2}\delta^4 u_j
\nonumber\\&&{}+\frac{2\gamma^3}{3h^4}\delta^6 u_j
-\frac{\gamma^3\alpha}{90h^2}\delta^6 u_j
\nonumber\\&&{}+\frac{\gamma^2 \alpha}{12h}\left( 2u_j\delta^3\mu
u_j+\delta^2 u_j \delta^3\mu u_j +\delta^4 u_j \delta \mu u_j
\right)\nonumber\\&&{}-\frac{\gamma^3\alpha}{480h}\left(16\, u_j\delta^5\mu
u_j +30\,\delta^4 u_j\delta^3\mu u_j +40\,\delta^2u_j\delta^3\mu u_j
\right.\nonumber\\&&\quad{}\left.
{}+40\,\delta^4 u_j \delta\mu u_j
+28\,\delta^2 u_j \delta^5\mu u_j +14\,\delta^6 u_j \delta\mu u_j
\right.\nonumber\\&&\quad{}\left.
{}+7\delta^4 u_j \delta^5\mu u_j +7\delta^6 u_j\delta^3\mu u_j\right)
+{\cal O}\left(\gamma^4,\alpha^2\right)\,.
\label{E_hol8dg2r1}
\end{eqnarray}
This discretisation forms a 7~point stencil approximation,
involving~$u_j$, $u_{j\pm1}$, $u_{j\pm2}$~and~$u_{j\pm3}$.
The first two lines of~(\ref{E_hol8dg2r1}), when evaluated at
$\gamma=1$, give a 4th~order centered difference approximation to the
hyperdiffusion term, a 6th~order centered difference approximation to
the linear growth term, and a 2nd~order centered difference
approximation to the nonlinear advection term.
The third and remaining lines account for higher order subgrid scale
dynamics of the nonlinearity and its inter-element coupling to generate
a 4th~order centered difference approximation to the nonlinearity~$uu_x$.

\paragraph{The $\Ord{\gamma^5,\alpha^2}$ holistic discretisation} is
\begin{eqnarray}
\dot{u_j}&=&-\frac{\gamma \alpha}{h^2}\delta^2 u_j
-\frac{\gamma \alpha}{h}u_j\delta\mu u_j
-\frac{4\gamma^2}{h^4}\delta^4 u_j
+\frac{\gamma^2 \alpha}{12h^2}\delta^4 u_j
\nonumber\\&&{}+\frac{2\gamma^3}{3h^4}\delta^6 u_j
-\frac{\gamma^3\alpha}{90h^2}\delta^6 u_j
\nonumber\\&&{}-\frac{7\gamma^4}{60h^4}\delta^8 u_j
+\frac{\gamma^4\alpha}{560h^2}\delta^8 u_j
\nonumber\\&&{}+\frac{\gamma^2 \alpha}{12h}\left(
2u_j\delta^3\mu u_j+\delta^2 u_j
\delta^3\mu u_j +\delta^4 u_j \delta \mu u_j \right)
\nonumber\\&&{}-\frac{\gamma^3\alpha}{480h}\left(
16\,u_j\delta^5\mu u_j +30\,\delta^4 u_j\delta^3\mu u_j
+40\,\delta^2u_j\delta^3\mu u_j
\right.\nonumber\\&&\quad{}\left.
{}+40\,\delta^4 u_j \delta\mu u_j
+28\,\delta^2 u_j
\delta^5\mu u_j +14\,\delta^6 u_j \delta\mu u_j
\right.\nonumber\\&&\quad{}\left.
{}+7\delta^4 u_j \delta^5\mu u_j
+7\delta^6 u_j\delta^3\mu u_j\right)
\nonumber\\&&{}+\frac{\gamma^4\alpha}{60480h} \left(
432\,u_j\delta^7\mu u_j +3528\,\delta^2 u_j \delta^5\mu u_j
+1507\,\delta^2u_j\delta^7 \mu u_j
\right.
\nonumber\\&&\quad{}\left.
{}+3780\,\delta^4u_j\delta^3\mu u_j
+3951\,\delta^4u_j\delta^5\mu u_j +984\,\delta^4 u_j \delta^7\mu u_j
\right.
\nonumber\\&&\quad{}\left.
{}+1764\,\delta^6u_j\delta\mu u_j
+3419\,\delta^6 u_j \delta^3\mu u_j
+1414\,\delta^6u_j\delta^5\mu u_j
\right.
\nonumber\\&&\quad{}\left.
{}+ 164\,\delta^6 u_j \delta^7\mu u_j +523\,\delta^8 u_j \delta\mu u_j
+656\,\delta^8 u_j \delta^3\mu u_j
\right.
\nonumber\\&&\quad{}\left.
{}+164\,\delta^8 u_j \delta^5\mu u_j
\right)+{\cal O}\left(\gamma^5,\alpha^2\right)\,.
\label{E_hol8dg3r1}
\end{eqnarray}
This forms a 9~point stencil approximation, involving only~$u_j$,
$u_{j\pm1}$, $u_{j\pm2}$, $u_{j\pm3}$~and~$u_{j\pm4}$.
The first two lines of~(\ref{E_hol8dg3r1}) when evaluated at $\gamma=1$
give a 6th~order centered difference approximation for the
hyperdiffusion term, an 8th~order centered difference approximation for
the linear growth term and a 2nd~order centered difference
approximation for the nonlinear advection term.
The third and remaining lines provide modifications to model the
nonlinear~$uu_x$ to 6th~order through resolving subgrid scale
dynamics.

We do not code these discretisations manually.
Instead, the computer algebra program at~\cite{Roberts02a} is
used with the \textsc{unix} editor \texttt{sed} to automatically write
the discretisation in a form suitable to be input to \textsc{Matlab}
for numerical exploration.

\paragraph{Compare to conventional centered difference models.}
Traditional direct finite differences generate the following
approximations to the \KS\ \pde~(\ref{EpdeKS}):
\begin{itemize}
    \item  5~point,
\begin{eqnarray}
\dot{u_j}&=&-\frac{\alpha}{h}u_j\delta\mu u_j
-\frac{\alpha}{h^2}\delta^2 u_j
-\frac{4}{h^4}\delta^4 u_j
+{\cal O}\left(h^2\right)\,;
\label{E_con2_ks}
\end{eqnarray}

    \item  7~point,
\begin{eqnarray}
\dot{u_j}&=&-\frac{\alpha}{h}\left(u_j\delta\mu u_j -\frac16\,u_j\delta^3\mu u_j\right)
-\frac{\alpha}{h^2}\left(\delta^2 u_j-\frac{1}{12}\,\delta^4 u_j\right)
\nonumber\\&&{}-\frac{4}{h^4}\left(\delta^4 u_j-\frac16\,\delta^6 u_j\right)+{\cal O}\left(h^4\right)\,;
\label{E_con4_ks}
\end{eqnarray}

    \item  9~point,
\begin{eqnarray}
\dot{u_j}&=&-\frac{\alpha}{h}\left(u_j\delta\mu u_j -\frac16\,u_j\delta^3\mu u_j
+\frac{1}{30}\,u_j\delta^5\mu u_j\right)
\nonumber\\&&{}
-\frac{\alpha}{h^2}\left(\delta^2 u_j-\frac{1}{12}\,\delta^4 u_j +\frac{1}{90}\,\delta^6 u_j\right)
\nonumber\\&&{}-\frac{4}{h^4}\left(\delta^4 u_j-\frac16\,\delta^6 u_j
+\frac{7}{240}\,\delta^6 u_j\right)+{\cal O}\left(h^6\right)\,.
\label{E_con6_ks}
\end{eqnarray}

\end{itemize}

Consider the different view of the errors for the discretisations:
the centered difference
approximations~(\ref{E_con2_ks}--\ref{E_con6_ks}) are justified by
consistency as grid size $h\to 0$\,; whereas the holistic
discretisations~(\ref{E_hol8dg1r1}--\ref{E_hol8dg3r1}) are supported by
centre manifold theory at finite grid size~$h$.
The errors in the centre manifold approach are due to the truncation of
dependence in the inter-element coupling parameter~$\gamma$ and the
nonlinearity parameter~$\alpha$.
However, as argued by Roberts~\cite{Roberts00a} for linear systems and
as demonstrated in \S\ref{S_ks_epde}, the particular choice of the
\ibc~(\ref{EbcsddKS}--\ref{EsbcidKS}) ensures that the holistic
discretisations \emph{are also consistent as $h\to 0$} with the \KS\
\pde~(\ref{EpdeKS}).

\subsection{Illustration of subgrid field enhances our view}
\label{S_ks_sg}

Recall the collection of subgrid fields (\ref{E_ks_cm_ab}) over the
physical domain form a state on the centre manifold.
Here we plot some example subgrid fields for various holistic models.
In particular, we examine subgrid fields of the holistic models of
steady states of the \KS\ \pde~(\ref{EpdeKS}) at nonlinear parameter
$\alpha=20$ and $\alpha=50\,$.
This is intended to reinforce the link between the abstract centre
manifold description of the dynamics and the physical subgrid fields
for the low order holistic models.
We compare the fields with the Lagrangian interpolation that underlies
traditional finite differences.
Recall that the key methodology difference is that the subgrid fields
of the holistic models are constructed by actual solutions of the \KS\
\pde, see~\S\ref{S_ks_sym}.

We restrict attention to odd symmetric solutions that are
$2\pi$-periodic.
This is done to compare with the numerical investigations of
Jolly~\cite{Jolly90} which we consider in more detail in
Sections \ref{chap_ks_ss}~and~\ref{chap_ks_td}.
Set a grid of 8~equi-spaced elements on the interval $[0,\pi]$.
The subgrid fields are plotted for approximations to the steady states
of the \KS\ equation~(\ref{EpdeKS}) with these periodic boundary
conditions, computed using holistic discretisations at $\alpha=20$ and
$\alpha=50\,$.

\begin{figure}
\centering
\includegraphics[width=0.8\textwidth]{kuraf_a20_1}
\caption{Subgrid field (green curve) of the holistic
model (\ref{E_hol8dg1r1}) and a
Lagrangian interpolant (magenta curve) constructed
through a 2nd~order centered difference approximation for a steady state of the \KS\ equation at
$\alpha=20$, with 8~elements on $[0,\pi]$.
An accurate solution is also plotted in blue.}
\label{kuraf_a20_1}
\end{figure}

Figure~\ref{kuraf_a20_1} displays an accurate solution
(blue curve) of the \KS\ \pde\ to compare with the subgrid field
(green curve) of the 5~point stencil
$\Ord{\gamma^3,\alpha^2}$ holistic approximation~(\ref{E_hol8dg1r1})
(green discs), and the Lagrangian interpolant
(magenta curve) constructed through a 2nd~order
centered difference approximation (magenta discs),
for a steady state at $\alpha=20\,$.
Observe the collection of subgrid fields forms the field~$u$ which is a
state on the centre manifold.
The subgrid field of the holistic model more accurately represents the
steady state of the \KS\ equation at $\alpha=20\,$, on this coarse
grid.

\begin{figure}
\centering
\includegraphics[width=0.8\textwidth]{kuraf_a20_2}
\caption{Subgrid fields of the holistic models with errors
$\Ord{\gamma^3,\alpha^2}$~(\ref{E_hol8dg1r1}) (green),
$\Ord{\gamma^4,\alpha^2}$~(\ref{E_hol8dg2r1})
(olive green) and $\Ord{\gamma^5,\alpha^2}$~(\ref{E_hol8dg3r1})
(cyan), for a steady state of the \KS\ equation at
$\alpha=20$, with 8~elements on $[0,\pi]$.
An accurate solution is also plotted in blue.}
\label{kuraf_a20_2}
\end{figure}

Higher order holistic models improve the accuracy and continuity of the
subgrid field.  Figure~\ref{kuraf_a20_2} displays the subgrid fields of
three holistic models for the same steady state of the \KS\ \pde\
depicted in Figure~\ref{kuraf_a20_1} for $\alpha=20$.  The
$\Ord{\gamma^3,\alpha^2}$ holistic model~(\ref{E_hol8dg1r1}) (green) is
the least accurate and has the largest jump at element boundaries.  The
$\Ord{\gamma^4,\alpha^2}$~(\ref{E_hol8dg2r1}) model (olive green)
displays improvement over the holistic $\Ord{\gamma^3,\alpha^2}$
approximation.  The $\Ord{\gamma^5,\alpha^2}$~(\ref{E_hol8dg3r1}) model
(cyan) is the most accurate, being almost indistinguishable from the
correct curve.

\begin{figure}
\centering 
\includegraphics[width=0.8\textwidth]{kuraf_a50_2}
\caption{Subgrid fields of the holistic models with errors
$\Ord{\gamma^3,\alpha^2}$~(\ref{E_hol8dg1r1}) (green),
$\Ord{\gamma^4,\alpha^2}$~(\ref{E_hol8dg2r1}) (olive green) 
and $\Ord{\gamma^5,\alpha^2}$~(\ref{E_hol8dg3r1})~(cyan), for a steady state of the \KS\ equation at
$\alpha=50$, with 8~elements on $[0,\pi]$.
An accurate solution is also plotted in blue.}
\label{kuraf_a50_1}
\end{figure}
Figure~\ref{kuraf_a50_1} shows a steady state of the \KS{} \pde\ at
$\alpha=50\,$.
The accurate field is symmetric (blue curve).
For this value of the nonlinearity there is no steady state solution
for centered difference approximations of either 2nd~(\ref{E_con2_ks}),
4th~(\ref{E_con4_ks}) or 6th~order~(\ref{E_con6_ks}) on this coarse
grid of 8~elements on $[0,\pi]$.
However, the 5~point stencil holistic approximation with errors
$\Ord{\gamma^3,\alpha^2}$~(\ref{E_hol8dg1r1}) (green)
models this steady state of the \KS\ equation even for such a large
value of the nonlinearity on this coarse grid.
This $\Ord{\gamma^3,\alpha^2}$ holistic solution has significant jumps
across the subgrid field at element boundaries; moreover, the subgrid
field is not symmetric and is most inaccurate near the centre of the
spatial domain considered here.
The 7~point stencil holistic approximation with errors
$\Ord{\gamma^4,\alpha^2}$~(\ref{E_hol8dg2r1})
(olive green) is more accurate with smaller jumps between neighbouring
the subgrid fields, but is also not symmetric.
The 9~point stencil holistic approximation with errors
$\Ord{\gamma^5,\alpha^2}$~(\ref{E_hol8dg3r1}) (cyan) is
the most accurate of the holistic models illustrated here; it is
symmetric and the jumps between neighbouring subgrid fields are almost
indiscernible.

These illustrations of the subgrid fields of steady states of the \KS{}
equation at $\alpha=20$ and $\alpha=50$ indicate the holistic models
perform well even at such large values of a supposedly small
parameter.  The performance of the holistic models are explored
further in Section~\ref{chap_ks_ss} for steady states and
Section~\ref{chap_ks_td} for time dependent phenomena.

\subsection{The holistic discretisations are consistent}
\label{S_ks_epde}
Holistic models constructed by implementing the
\ibc~(\ref{EbcsddKS}--\ref{EsbcidKS}) have dual
justification~\cite{Roberts00a}: they are supported by centre manifold
theory for small enough $\alpha$~and~$\gamma$; as well as being
justified by their consistency as the grid size~$h\to0$.  We explore
consistency as a well established feature of numerical
analysis.\footnote{But note that high order consistency is not a
primary goal of this holistic approach, since we aim to develop and
support models for finite element size~$h$.}

Here we examine the equivalent \pde{}s for the holistic
discretisations~(\ref{E_hol8dg1r1}--\ref{E_hol8dg3r1}) evaluated at
$\gamma=1$, and the centered difference
approximations~(\ref{E_con2_ks}--\ref{E_con6_ks}).
These equivalent \pde{}s establish the $\Ord{h^{2p-2}}$ consistency
with the \KS\ \pde\ for holistic models constructed with residuals
$\Ord{\gamma^{p+1}}$.

Roberts~\cite{Roberts00a} proved that using \ibc\ of the form
introduced in \S\ref{S_KSho} and retaining terms up to~$\gamma^p$ in
the holistic approximations results in approximations which are
consistent with the \emph{linear} terms of the \KS\ equation~(\ref{EpdeKS}) to
$\Ord{h^{2p-2}}$, provided $p\ge2$\,.
However, it appears that using the
\ibc~(\ref{EbcsddKS}--\ref{EsbcidKS}) also ensures $\Ord{h^{2p-2}}$
consistency for the nonlinear terms.
As yet no formal proof exists of this nonlinear consistency, but all
holistic models of the \KS\ equation, containing terms up to
$\gamma^7$~and~$\alpha^4$ and constructed
using~(\ref{EbcsddKS}--\ref{EsbcidKS}) are nonlinearly consistent
(although not all are recorded here).

Find the equivalent \pde s for the various discretisations by
expanding the discretisations in grid size~$h$ about a grid point~$x_j$.
That is, write
\begin{equation}
u_{j\pm m}=u_j\pm mh\D x{u_j}+m^2\frac{h^2}{2}\DD x{u_j}
+\sum^{\infty}_{k=3}(\pm m)^k\frac{h^k}{k!}\DG x{u_j}{k}\,,
\end{equation}
to whatever order in~$h$ is required.  Computer algebra performs the
tedious details.

\paragraph{The equivalent PDE for the 5~point holistic
discretisation} (\ref{E_hol8dg1r1}), which retain terms up to
$\gamma^2$, is
\begin{eqnarray}
\D tu &=& -\alpha\left(u\D xu +\DG xu2\right) -4\DG xu4
-\frac{2h^2}{3}\DG xu6+\frac{h^4}{20} \DG xu8
\nonumber\\&&{}
+\alpha h^4
\left(\frac{1}{48}\DG xu3\DG xu2 +\frac{1}{48}\DG xu4 \D xu
+\frac{1}{30}u\DG xu5+\frac{1}{90}\DG xu6\right)
\nonumber\\&&{}
+{\cal O}(h^6)\,.\label{E_epde_holdg1r1}
\end{eqnarray}
The equivalent \pde\ for the 5~point centered difference
approximation~(\ref{E_con2_ks}) is
\begin{eqnarray}
\D tu &=& -\alpha\left(u\D xu +\DG xu2\right) -4\DG xu2
\nonumber\\&&{}-h^2\left(\alpha\frac16\DG xu3 \D xu+\alpha\frac{1}{12}\DG xu4+ \frac23\DG xu6\right)
\nonumber\\&&{}
+{\cal O}(h^4)\,.
\label{E_epde_con2_ks}
\end{eqnarray}
Observe that both equivalent \pde s
(\ref{E_epde_holdg1r1}--\ref{E_epde_con2_ks}) are $\Ord{h^2}$ accurate.
The coefficients of the error terms are different in both of the these
equivalent \pde s, with those of~(\ref{E_hol8dg1r1}) having fewer
error terms.

\paragraph{The equivalent PDE for the 7~point holistic
discretisation} (\ref{E_hol8dg2r1}), which
retains terms up to~$\gamma^3$, is
\begin{eqnarray}
\D tu &=& -\alpha\left(u\D xu +\DG xu2\right) -4\DG xu4
-\frac{7h^4}{60}\DG xu8+\frac{13h^6}{756} \DG{x}{u}{10}
\nonumber\\&&{}
-\alpha h^6
\left(\frac{17}{640}\DG xu4\DG xu3 +\frac{7}{384}\DG xu5 \DG xu2
+\frac{3}{320}\DG xu6 \D xu\right.
\nonumber\\&&{}\left.
+\frac{1}{140}u\DG xu7
+\frac{1}{560}\DG xu8\right)
\nonumber\\&&{}
+{\cal O}(h^8)\,.
\label{E_epde_holdg2r1}
\end{eqnarray}
Whereas the equivalent \pde\ for the 7~point centered difference
approximation~(\ref{E_con4_ks}) is
\begin{eqnarray}
\D tu &=& -\alpha\left(u\D xu +\DG xu2\right) -4\DG xu2
\nonumber\\&&{}-h^4\left(\alpha\frac{1}{30}\DG xu5 \D xu+\alpha\frac{1}{90}\DG xu6+ \frac{7}{60}\DG xu8\right)
\nonumber\\&&{}
+{\cal O}(h^6)\,.
\label{E_epde_con4_ks}
\end{eqnarray}
The two equivalent \pde{}s
(\ref{E_epde_holdg2r1}--\ref{E_epde_con4_ks}) are $\Ord{h^4}$
accurate, and again the holistic discretisation has fewer errors.

\paragraph{The equivalent PDEs for the 9~point holistic
discretisation} (\ref{E_hol8dg2r1}), which retain terms up to $\gamma^4$, 
is
\begin{eqnarray}
\D tu &=& -\alpha\left(u\D xu +\DG xu2\right) -4\DG xu4
-\frac{41h^6}{1890}\DG {x}{u}{10}+\frac{13h^8}{2700} \DG{x}{u}{12}
\nonumber\\&&{}
+\alpha h^8
\left(\frac{3433}{138240}\DG xu5\DG xu4
+\frac{5927}{322560}\DG xu6 \DG xu3
+\frac{499}{53760}\DG xu7 \DG xu2\right.
\nonumber\\&&{}\left.
+\frac{29}{8960}\DG xu8 \D xu
+\frac{1}{630}u\DG xu9
+\frac{1}{3150}\DG {x}{u}{10}
\right)
\nonumber\\&&{}
+{\cal O}(h^{10})\,.
\label{E_epde_holdg3r1}
\end{eqnarray}
The equivalent \pde\ for the 9pt~centered difference
approximation~(\ref{E_con6_ks}) is
\begin{eqnarray}
\D tu &=& -\alpha\left(u\D xu +\DG xu2\right) -4\DG xu2
\nonumber\\&&{}-h^6\left(\alpha\frac{1}{140}\DG xu7 \D xu+\alpha\frac{1}{560}\DG xu8+ \frac{41}{1890}\DG xu{10}\right)
\nonumber\\&&{}
+{\cal O}(h^8)\,.
\label{E_epde_con6_ks}
\end{eqnarray}
Again the the equivalent \pde s
(\ref{E_epde_holdg3r1}--\ref{E_epde_con6_ks}) are $\Ord{h^6}$ accurate,
with the holistic discretisation having fewer error terms.

Although there is no proof of nonlinear consistency in general, we
have demonstrated it here for these three holistic discretisations,
and have found nonlinear consistency for all models investigated.

\section{Holistic models accurately give steady states}
\label{chap_ks_ss}

The relevance of our holistic models is rigorously supported by centre
manifold theory for sufficiently small parameters
$\gamma$~and~$\alpha\,$.
However, the holistic models must be evaluated at coupling parameter
$\gamma=1$ to model the dynamics of the \KS\ equation.
The important question: Does evaluating the holistic models at
$\gamma=1$ provide useful and accurate numerical models?
Numerical experiments detailed in this and the next section provide
strong support that it does.

In this section we explore the accuracy of the holistic models by
constructing and comparing bifurcation diagrams of the various holistic
discretisations to conventional explicit centered difference
approximations and to the bifurcation diagrams presented by Jolly
et~al.~\cite{Jolly90} for various traditional Galerkin and nonlinear
Galerkin approximations.

We restrict exploration to solutions that are both $2\pi$~periodic and
odd: thus
\begin{equation}
u(x,t)=u(x+2\pi,t) \quad \mbox{and} \quad u(x,t)=-u(2\pi-x,t)\,.
\label{E_odd_bc}
\end{equation}
We also restrict the nonlinearity parameter to the range $0\le \alpha
\le 70$\,.
These restrictions are to compare our results to those of Jolly
et~al.~\cite{Jolly90} for approximate inertial manifold methods.
For this range of nonlinearity~$\alpha$ the trivial solution $u=0$
undergoes pitchfork bifurcations at $\alpha=4,16,36,64$ leading to the
unimodal, bimodal, trimodal and quadrimodal branches respectively, see
the bifurcation diagram Figure~\ref{bif6ord48}.

\begin{figure}
\centering
\includegraphics[width=0.95\textwidth]{con6ord48}
\caption{Accurate bifurcation diagram $0\le\alpha\le70$ for the
\KS\ equation, using a 6th~order centered difference
approximation with 48~points on the interval~$[0,\pi]$. A signed
$L^2$~norm is plotted against~$\alpha$}
\label{bif6ord48}
\end{figure}

Such bifurcation diagrams usefully summarise qualitative and
quantitative information for a large range of the nonlinearity
parameter~$\alpha$.
We use the software package \textsc{xppaut}~\cite{xppaut01}, which
incorporates the continuation software \textsc{auto}~\cite{auto01}, to
calculate the bifurcation information.
The information is then filtered through a function written in
\textsc{matlab} to draw the bifurcation diagram.
The input to \textsc{xppaut} is a text \verb|.ode| file describing the
set of \textsc{ode}s.
Because the holistic models contain a large number of terms the
\verb|.ode| files are generated automatically using \textsc{reduce} and
\textsc{matlab} which also incorporates the odd periodic
requirement~(\ref{E_odd_bc}), see \cite{MacKenzie05} for more details.

\subsection{Reference accurate steady states}
\label{S_numksss}

Here we introduce accurate solutions for the steady states of the \KS{}
equation~(\ref{EpdeKS}) over the range $0\le \alpha \le70$\, as
summarised in the bifurcation diagram of Figure~\ref{bif6ord48}.
Accurate solutions are produced by a 6th~order accurate centered
difference approximation~(\ref{E_con6_ks}) with 48~grid points on the
spatial interval $[0,\pi]$.
These provide the reference for the approximations on coarse grids,
and serve to also introduce the conventions we adopt in bifurcation
diagrams.

For all the bifurcation diagrams a signed solution norm is plotted
against the nonlinearity parameter~$\alpha$.
This is different to the convention adopted by Jolly
et~al.~\cite{Jolly90} but empowers us to investigate more detail by
showing positive and negative branches---stability differs along these
branches.
For example, see in Figure~\ref{bif6ord48} that the negative bimodal
branch is stable for $16.140<\alpha<22.556$, whereas the positive
bimodal branch is unstable.
The solution norm is signed corresponding to the sign of the the first
grid value, $u_1=u(x_1)\,$.
The blue curves are branches of stable fixed points and the red curves
are branches of unstable fixed points.
The open squares denote pitchfork bifurcations and the black squares
denote Hopf bifurcations.

The labeling scheme used in Figure~\ref{bif6ord48} follows that of
Jolly et~al.~\cite{Jolly90} and Scovel~\cite{Scovel88} with the
addition of a plus or minus sign depending upon the sign of~$u_1$.
For example, the secondary bifurcation on the negative bimodal branch
is labeled $R_2b_1-$ from the labeling scheme of Scovel with the
addition of the $-$~sign because it occurs on the negative branch.
Figure~\ref{bif6ord48} appears to show several discontinuities.
For example, the positive unimodal branch ends at approximately
$\alpha=12\,$.
This apparent discontinuity arises due to the convention adopted here
of taking the sign of~$u_1$ to sign the norm: actually there is a
continuous transformation as the positive unimodal branch and the
negative unimodal branch transform into the negative bimodal branch.
It is straightforward to sign the branch near the trivial solution, but
away from the trivial solution the distinction between positive and
negative may be ambiguous and occasionally leads to jumps in the
bifurcation diagram.

\begin{figure}
\centering
\includegraphics[width=0.95\textwidth]{ks_sol_stable}
\caption{Some examples of the stable equilibria of the \KS\ equation.
Dark blue curves are solutions along the negative unimodal and bimodal
branches.  Light blue curves are stable solutions along the negative
trimodal branch.}
\label{ks_sol_stable}
\end{figure}

For later comparison see in Figure~\ref{ks_sol_stable} some of the 
stable equilibria of the \KS\ equation in the regime of interest,
$0\le\alpha\le70\,$.
Figures~\ref{ks_sol_stable}a,b,c show solutions on the negative
unimodal branch at $\alpha=1,5,10$ respectively.
Figures~\ref{ks_sol_stable}d,e,f show solutions on the negative
bimodal branch at $\alpha=20,30,40$ respectively.
The dark blue curves in Figures~\ref{ks_sol_stable}g,h,i show
solutions on the negative bimodal branch and the light blue curves are
solutions on the negative trimodal branch at $\alpha=50,55,60$
respectively.

\subsection{Holistic models are accurate on coarse grids}
\label{S_ks_ss_coarse}

We begin investigating the performance of the holistic models by
considering the ${\cal O}(\gamma^5,\alpha^2)$ holistic
model~(\ref{E_hol8dg3r1}) (9~point stencil, $\Ord{h^6}$ consistent).
We investigate its reproduction of the steady states of the
\KS\ system using coarse grids on the interval $[0,\pi]$.

\begin{figure}
\centering
\includegraphics[width=0.95\textwidth]{ks_acc1}
\caption{Some accurate solutions plotted with holistic and centered
difference approximations on coarse grids.  Blue curves are accurate
solutions, green curves are the holistic approximation with 
\textsc{ibc}s~(\ref{EbcsddKS}--\ref{EsbcidKS}) with errors
$\Ord{\gamma^5,\alpha^2}$ on 8~elements.  Magenta curves are a
6th~order centered difference approximation with 8~grid points.}
\label{ks_acc1}
\end{figure}

Figure~\ref{ks_acc1} shows accurate solutions of the \KS\
equation~(\ref{EpdeKS}) with odd boundary conditions~(\ref{E_odd_bc})
in blue.
The holistic model with errors ${\cal O}(\gamma^5,\alpha^2)$ and
8~elements is shown in green: Figure~\ref{ks_acc1}g,h,i, the bottom
row, shows that the holistic model with errors ${\cal
O}(\gamma^5,\alpha^2)$ gives at large nonlinearity the stable bimodal
and trimodal solutions, $\alpha=50$ and $\alpha=55$\,, and the stable
bimodal solution at $\alpha=60$\,.
Magenta curves are solutions of the 6th~order centered difference
approximation~(\ref{E_con6_ks}) with 8~grid points---it has equal
stencil width to the holistic model.
The 6th~order centered difference approximation does not give any
stable solutions for $\alpha\geq50$\,.
The holistic model provides reasonable solutions where comparable
traditional methods do not.

\subsubsection{Bifurcation diagrams show success}
\label{S_ks_ss_cg}

Now turn to the bifurcation diagram to obtain a more comprehensive view.
We see the holistic model has good bifurcation diagrams on a coarse
grid of 8~elements and even with just 6~elements.

\begin{figure}
\centering
\includegraphics[width=0.95\textwidth]{hol_high_8bif_1}
\caption{Bifurcation diagrams for coarse grid approximations with
8~elements on $[0,\pi]$ for (a)~holistic model ${\cal
O}(\gamma^5,\alpha^2)$, (b)~centered difference
6th~order.}
\label{hol_high_8bif_1}
\end{figure}

Figure~\ref{hol_high_8bif_1} shows a side by side comparison of the
holistic model with  errors
$\Ord{\gamma^5,\alpha^2}$ with 8~elements on $[0,\pi]$ and the
6th~order centered difference approximation with 8~grid points on
$[0,\pi]$.
These approximations are both 9~point stencil approximations.
The accurate bifurcation diagram is also plotted in grey but without
any stability information.
The signed $L_2$~norms for the bifurcation diagrams on the coarse grid
of 8~elements are adjusted by a factor of~$\sqrt{6}$ to allow
comparison to the accurate bifurcation diagram constructed with 48~grid
points on $[0,\pi]$.
Throughout this paper when comparing bifurcation diagrams of
different grid resolutions, the signed $L_2$ norms are adjusted this
way to provide a consistent reference.
Figure~\ref{hol_high_8bif_1}a shows the $\Ord{\gamma^5,\alpha^2}$
holistic model gives good agreement with the accurate bifurcation
diagram for $\alpha<40$ and qualitatively reproduces most of the
bifurcation picture for $40<\alpha<70\,$.
The $\Ord{\gamma^5,\alpha^2}$ holistic model does not detect the
bifurcation points $R_3t_2\pm$ on this coarse grid and the bifurcation
points $R_3t_1\pm$ are incorrectly identified as fold points.
However, the $\Ord{\gamma^5,\alpha^2}$ holistic model finds all of the
other bifurcation points in this range of~$\alpha$.
Figure~\ref{hol_high_8bif_1}b shows the 6th~order centered difference
approximation gives good agreement with the accurate bifurcation
diagram only for $\alpha<20$ and qualitatively reproduces the
bifurcation diagram for $20<\alpha<40\,$.
The 6th~order centered difference approximation performs poorly for
$\alpha>40\,$.
Table~\ref{tablebifs} lists the values of~$\alpha$ at which the
bifurcation points occur and confirms the $\Ord{\gamma^5,\alpha^2}$
holistic model performs more accurately than the 6th~order centered
difference approximation on this coarse grid of 8~elements.

\begin{table}
\caption{$\alpha$ values at which bifurcation points occur for the
various coarse grid approximations;
$^*$~denotes bifurcation point identified as fold point.}
\label{tablebifs}
\begin{center}\footnotesize
\begin{tabular}{lllllllll}
\hline
Approximation&$R_2b_1$&$R_2b_2$&$R_2b_3$&$R_2b_4$&$R_3t_1$&$R_3t_2$&$R_4b_1$&
$R_4q_1$\\ 
\hline
\multicolumn{4}{l}{Accurate $48$\,pts}\\ 
$6\mbox{th}$ order
&16.14&22.56&52.89&63.74&36.23&50.91&64.56&64.28\\ 
\hline
\multicolumn{4}{l}{Holistic~8 elements}
\\ 
$\Ord{\gamma^3,\alpha^2}$&14.64&20.36&39.34&44.96&$29.28^*$
&---&45.28&44.87\\ 
$\Ord{\gamma^3,\alpha^3}$&14.65&20.52&39.66&45.16&$29.33^*$
&---&45.47&44.96\\ 
$\Ord{\gamma^3,\alpha^4}$&14.65&20.53&39.72&45.21&$29.33^*$
&---&45.51&44.97\\ 
$\Ord{\gamma^4,\alpha^2}$&16.00&22.56&48.62&57.38&$34.73^*$&
---&57.89&57.49\\ 
$\Ord{\gamma^4,\alpha^3}$&16.00&22.56&48.25&56.84&$34.73^*$&
---&57.45&57.28\\ 
$\Ord{\gamma^4,\alpha^4}$&16.00&22.57&48.10&56.63&$34.73^*$&
---&57.30&57.21\\ 
$\Ord{\gamma^5,\alpha^2}$&16.13&22.72&51.54&61.54&$35.89^*$&
---&62.20&61.78\\ 
$\Ord{\gamma^5,\alpha^3}$&16.13&22.73&51.53&61.37&$35.91^*$&
---&62.04&61.70\\ 
$\Ord{\gamma^5,\alpha^4}$&16.13&22.73&51.60&61.38&$35.91^*$&
---&62.02&61.69\\ 
\hline
\multicolumn{4}{l}{Centered 8~pts}\\ 
$2\mbox{nd}$ order&15.30&19.81&---&---&---&---&---&---\\ 
$4\mbox{th}$ order&16.02&21.55&---&---&$35.94^*$&---&---&---\\ 
$6\mbox{th}$ order&16.12&21.99&---&---&$35.83^*$&---&---&---\\ 
\hline
\multicolumn{4}{l}{Holistic~12 elements}\\ 
$\Ord{\gamma^3,\alpha^2}$&15.45&21.67&45.96&53.94&32.86
&45.98&54.49&54.17\\ 
$\Ord{\gamma^3,\alpha^3}$&15.45&21.69&46.05&54.00&32.87
&46.33&54.55&54.20\\ 
$\Ord{\gamma^4,\alpha^2}$&16.11&22.62&51.93&62.10&35.90
&50.92&62.83&62.52\\ 
$\Ord{\gamma^4,\alpha^3}$&16.11&22.62&51.94&62.10&35.90
&50.94&62.83&62.52\\ 
\hline
\multicolumn{4}{l}{Centered 12~pts}\\ 
$2\mbox{nd}$ order&15.77&21.68&48.33&57.63&34.36&44.70&58.34&
58.36\\ 
$4\mbox{th}$ order&16.12&22.37&51.74&62.33&35.98&48.62&63.11&
62.98\\ 
\end{tabular}
\end{center}
\end{table}

\begin{figure}
\centering
\includegraphics[width=0.95\textwidth]{hol_high_6bif_1}
\caption{Bifurcation diagrams for coarse grid approximations with
6~elements on $[0,\pi]$ for (a) holistic model ${\cal
O}(\gamma^5,\alpha^2)$, (b) centered difference
6th~order.}
\label{hol_high_6bif_1}
\end{figure}
Figure~\ref{hol_high_6bif_1} is a side by side comparison of the same
$\Ord{\gamma^5,\alpha^2}$ holistic model to the
6th~order centered difference approximation, on an even coarser grid of
just 6~elements.
The superior performance of the holistic model is again evident.
We conjecture that the superior performance of the holistic
discretisation is due to its systematic modelling of the subgrid scale
processes.
These bifurcation diagrams, Figures~\ref{hol_high_8bif_1}
and~\ref{hol_high_6bif_1}, give excellent support to the holistic
approach to generating approximations for the \KS\ equation.

We also investigate various holistic models for the \KS\ \pde\ by
comparing bifurcation diagrams of holistic models of higher orders.
We examine bifurcation diagrams for holistic models
with errors ${\cal O}(\gamma^p,\alpha^q)$, for $3\le p \le 5$ and
$2\le q \le 4$\,, and find that retaining terms of higher order in coupling
parameter~$\gamma$, corresponding to wider stencil approximations,
gives much greater improvement in accuracy than retaining terms of
higher order in the nonlinearity parameter~$\alpha$.

\begin{figure}
\centering
\includegraphics[width=0.95\textwidth]{hol8d_bif}
\caption{Bifurcation diagrams
for the holistic models with 8~elements on the interval
$[0,\pi]$ up to and including the ${\cal
O}(\gamma^5,\alpha^4)$ holistic model.}
\label{hol8d_bif}
\end{figure}

Figure~\ref{hol8d_bif} shows the bifurcation diagrams for the holistic
models up to and including the ${\cal
O}(\gamma^5,\alpha^4)$ holistic model.
Surveying across the columns of Figure~\ref{hol8d_bif} see the
bifurcation diagrams for holistic models of increasing order of
coupling parameter~$\gamma$, corresponding to approximations of
increasing stencil width.
For example, Figure~\ref{hol8d_bif}a,b,c shows the bifurcation
diagrams for the holistic models~(\ref{E_hol8dg1r1}),
(\ref{E_hol8dg2r1}) and~(\ref{E_hol8dg3r1}) respectively.
Surveying down the rows of Figure~\ref{hol8d_bif}  see the
bifurcation diagrams for increasing orders of the nonlinearity
parameter~$\alpha$.
Figure~\ref{hol8d_bif} illustrates the improvement in accuracy of the
higher order holistic models.
Note first the dramatic improvement in accuracy gained by
moving from left to right across Figure~\ref{hol8d_bif}, corresponding
to approximations of higher orders in the coupling parameter~$\gamma$.

Second, see that less improvement is gained by moving from top to
bottom of Figure~\ref{hol8d_bif}, corresponding to approximations of
higher order in the nonlinearity parameter~$\alpha$.
There are some peculiarities about this series of bifurcation pictures
for holistic models of increasing order in $\alpha$.
For the 5~point stencil approximations displayed in the first column of
Figures~\ref{hol8d_bif}, higher orders in~$\alpha$ appear to gain some
improvement.
In particular Figures~\ref{hol8d_bif}d,g show the ${\cal
O}(\gamma^3,\alpha^3)$ and ${\cal O}(\gamma^3,\alpha^4)$ holistic
models reproduce the unstable trimodal branches that were missing from
the ${\cal O}(\gamma^3,\alpha^2)$ bifurcation diagram shown in
Figure~\ref{hol8d_bif}a.
However, for the 7~point stencil approximations displayed in the second
column of Figure~\ref{hol8d_bif}, holistic models of higher orders
in~$\alpha$ lose some features of the \KS\ system.
The correct behaviour of the unstable trimodal and quadrimodal branches
is reproduced for the ${\cal O}(\gamma^4,\alpha^2)$ model shown in
Figure~\ref{hol8d_bif}b, but not reproduced for the higher order
${\cal O}(\gamma^4,\alpha^3)$ and ${\cal O}(\gamma^4,\alpha^4)$ models
shown in Figures~\ref{hol8d_bif}e,h respectively.
For the 9~point stencil approximations, displayed in the third column
of Figures~\ref{hol8d_bif}, the $\Ord{\gamma^5,\alpha^2}$ holistic
model shown in Figure~\ref{hol8d_bif}c, reproduces the unstable
trimodal branch whereas the higher order $\Ord{\gamma^5,\alpha^3}$
model shown in Figure~\ref{hol8d_bif}f, does not reproduce the
unstable trimodal branch.
These peculiarities suggest that while we have observed excellent
performance of the holistic models constructed with the non-local~\ibc\
on coarse grids, it may be possible that modifications could be made to
the non-local \ibc\ such that higher order approximations in the
nonlinear parameter are improved.
Exploration of possible such modifications are left for further
research.

\subsubsection{Holistic models outperform centered differences}

In \S\ref{S_ks_ss_cg} we saw that the performance of the
$\Ord{\gamma^5,\alpha^2}$ holistic model~(\ref{E_hol8dg3r1})
constructed with non-local \ibc\ was far superior to the explicit
6th~order centered difference approximation (\ref{E_con6_ks}).  To
complete the comparison of holistic models to explicit centered
difference schemes, we compare the
$\Ord{\gamma^3,\alpha^2}$~(\ref{E_hol8dg1r1}) and
$\Ord{\gamma^4,\alpha^2}$~(\ref{E_hol8dg2r1}) holistic models to the
2nd~order~(\ref{E_con2_ks}) and 4th~order~(\ref{E_con4_ks}) centered
difference approximations respectively; these are 5~point and 7~point
discretisations respectively.

\begin{figure}
\centering
\includegraphics[width=0.95\textwidth]{con_ks_8bif_all}
\caption{Bifurcation diagrams for (a)~$\Ord{\gamma^3,\alpha^2}$
holistic model, (b)~2nd~order centered difference,
(c)~$\Ord{\gamma^4,\alpha^2}$ holistic model and (d)~4th~order centered
difference all with 8~elements on the interval~$[0,\pi]$}
\label{con_ks_8bif_all}
\end{figure}

The first row of Figure~\ref{con_ks_8bif_all} is a side by side
comparison of the $\Ord{\gamma^3,\alpha^2}$ holistic model and the
2nd~order centered difference approximation with 8~elements on
$[0,\pi]$.  The second row of Figure~\ref{con_ks_8bif_all} is a side by
side comparison of the $\Ord{\gamma^4,\alpha^2}$ holistic model and the
4th~order centered difference approximation on the same coarse grid.
The accurate bifurcation diagram is plotted in grey without any
stability information.

Although comparing Figures~\ref{con_ks_8bif_all}b,d shows some
improvement is gained by taking higher order centered difference
approximations, this improvement is not as pronounced as for the
holistic models on this coarse grid as shown in
Figures~\ref{con_ks_8bif_all}a,c.
Both the 2nd~order and 4th~order centered difference approximations
fail to reproduce the correct behaviour of the unstable trimodal and
quadrimodal branches.
In contrast, even the 5~point stencil $\Ord{\gamma^3,\alpha^2}$
holistic approximation qualitatively reproduces the trimodal and
quadrimodal branches on the same coarse grid.
The values at which the bifurcation points occur are listed in
Table~\ref{tablebifs} and confirm these holistic models outperform the
centered difference approximations on this coarse grid of 8~elements
on~$[0,\pi]$.

\subsubsection{Grid refinement improves accuracy}

Since the equivalent \pde{}'s, (\ref{E_epde_holdg1r1}),
(\ref{E_epde_holdg2r1}) and~(\ref{E_epde_holdg3r1}), for our holistic
models are of $\Ord{h^2}$, $\Ord{h^4}$ and $\Ord{h^6}$ respectively,
grid refinement should result in improved accuracy.

\begin{figure}
\centering
\includegraphics[width=0.95\textwidth]{hol12d_bif1}
\caption{Bifurcation diagrams for the holistic models 
with 12~elements on the interval~$[0,\pi]$.  Compare with
Figure~\ref{hol8d_bif} with 8~elements.}
\label{hol12d_bif}
\end{figure}

Figure~\ref{hol12d_bif} shows the bifurcation diagrams of the holistic
models up to and including the $\Ord{\gamma^4,\alpha^3}$ model on a
finer grid of 12~elements on $[0,\pi]$.
Compare Figure~\ref{hol12d_bif} with Figure~\ref{hol8d_bif} to confirm
the improved accuracy for the holistic models on this refined grid.
Table~\ref{tablebifs} also shows the bifurcation points are more
accurately reproduced for the holistic models on this refined grid.

\begin{figure}
\centering
\includegraphics[width=0.95\textwidth]{hol12d_bif2}
\caption{Bifurcation diagrams for (a)~$\Ord{\gamma^4,\alpha^2}$
holistic model, (b)~4th~order centered difference approximations
with 12~elements on the interval~$[0,\pi]$}
\label{hol12d_bif2}
\end{figure}

Figure~\ref{hol12d_bif2} is a side by side comparison of the
bifurcation diagrams of the $\Ord{\gamma^4,\alpha^2}$ holistic model
and the 4th~order centered difference approximation~(\ref{E_con4_ks}).
The accurate bifurcation diagram is shown in grey.
See the $\Ord{\gamma^4,\alpha^2}$ holistic model is more accurate for
$0\le\alpha\le70$ but the improvement is not as pronounced as it is on
the coarser grid of 8~elements.
We suggest that this is because the major benefit to using the holistic
models comes from application on coarser grids where the subgrid scale
modelling is more significant.

\subsection{Comparison to Galerkin approximations}
\label{S_ks_ss_g}
Here we investigate the traditional Galerkin and non-linear Galerkin
approximations~\cite{Jolly90} for the \KS\ equation~(\ref{EpdeKS}) with
the periodic and odd conditions~(\ref{E_odd_bc}).
We find the holistic models compare well with the Galerkin methods.
While the Galerkin methods are of superior accuracy for solving the
\KS\ system~(\ref{EpdeKS}) with periodic boundary conditions, because
of their global nature they lack the flexibility of the local nature of
the holistic models.
Although not explored here, this local nature of the holistic models
empowers its use with physical boundary conditions~\cite{Roberts01b}
other than periodic.

Galerkin methods seek solutions in the form which is dominantly the
superposition of $m$~periodic, global modes:
\begin{equation}
u(x,t)=\sum^m_{k=1}b_k(t)\sin(kx)\,. 
\label{E_gal_u}
\end{equation}
\paragraph{The $m$-mode traditional Galerkin}
approximation~\cite{Jolly90} is
\begin{equation}
\frac{d b_k}{dt}\approx\left(-4k^4+\alpha k^2\right)b_k
-\alpha\beta^m_k\,,\quad 1\le k\le m\,,
\label{E_tr_gal_ks}
\end{equation}
where
\begin{equation}
\beta^m_k(b_1,\ldots,b_m)
=\frac 12 \sum^m_{j=1}jb_j\left[b_{k+j}+\mbox{sign}(k-j)b_{|k-j|}\right]\,.
\label{E_beta_ks}
\end{equation}

\paragraph{The $m$-mode first iterate nonlinear Galerkin}
approximation~\cite{Jolly90} is based upon the adiabatic
approximation~(\ref{E_nl_gal2_ks}) for higher wavenumber modes
$k=m+1:2m$, namely
\begin{equation}
\frac{d b_k}{dt}\approx\left(-4k^4+\alpha k^2\right)b_k
-\alpha\beta^{2m}_k(b_1,\ldots,b_m,\phi_{m+1}
,\ldots,\phi_{2m})
\,,
\label{E_nl_gal1_ks}
\end{equation}
for $1\le k\le m$\,, where
\begin{equation}
\phi_j=-\frac{\alpha}{4j^4}\beta^{2m}_j
\left(b_1,\ldots,b_m,0,\ldots,0\right)\,,
\label{E_nl_gal2_ks}
\end{equation}
for $m+1\le j \le 2m$ and $\beta^{2m}_j$ is given by~(\ref{E_beta_ks}).

Obtain higher order nonlinear Galerkin approximations~\cite{Roberts89}
through recognising time derivatives of these and even higher wave
number modes.
We do not explore these.

\begin{figure}
\centering
\includegraphics[width=0.95\textwidth]{tr_gal_bifs}
\caption{Bifurcation diagrams for (a)~3~mode, (b)~4~mode, (c)~6~mode
and (d)~8~mode traditional Galerkin approximations on~$[0,\pi]$.}
\label{nl_galt_bifs}
\end{figure}

Now examine the bifurcation diagrams of the two Galerkin
approximations~(\ref{E_gal_u}--\ref{E_nl_gal2_ks}) for
$0\le\alpha\le70$ and compare with the bifurcation diagrams of the
holistic models on coarse grids, presented in \S\ref{S_ks_ss_coarse}.
Figure~\ref{nl_galt_bifs} shows the Bifurcation diagrams for the
3~mode, 4~mode, 6~mode and 8~mode traditional Galerkin approximations
on~$[0,\pi]$.
See that at least 4~modes are needed to qualitatively reproduce the
behaviour of the stable bimodal branch.
Compare the $\Ord{\gamma^5,\alpha^2}$ holistic model with 6~elements
from Figure~\ref{hol_high_6bif_1}a, to the 6~mode traditional Galerkin
approximation and observe the $\Ord{\gamma^5,\alpha^2}$ holistic model
qualitatively models most steady state dynamics that are reproduced by
the 6~mode traditional Galerkin approximation.
Neither the $\Ord{\gamma^5,\alpha^2}$ holistic model nor the 6~mode
traditional Galerkin approximation qualitatively reproduce the correct
behaviour of the unstable quadrimodal branch.
Similarly the $\Ord{\gamma^5,\alpha^2}$ holistic model with 8~elements
from Figure~\ref{hol_high_8bif_1}a and the 8~mode traditional Galerkin
approximation qualitatively model most steady state dynamics.
However, the 8~mode traditional Galerkin approximation is more accurate.

\begin{figure}
\centering
\includegraphics[width=0.95\textwidth]{nl_gal_bifs}
\caption{Bifurcation diagrams for (a)~3~mode, (b)~4~mode, (c)~6~mode
and (d)~8~mode first iterate nonlinear Galerkin approximations on
$[0,\pi]$.}
\label{nl_gal_bifs}
\end{figure}

Figure~\ref{nl_gal_bifs} shows the bifurcation diagrams for the 3~mode,
4~mode, 6~mode and 8~mode first iterate nonlinear Galerkin
approximations~(\ref{E_nl_gal1_ks}) on~$[0,\pi]$.
See impressive accuracy for the low mode first iterate nonlinear
Galerkin approximations.
The 6~mode nonlinear Galerkin approximation reproduces all of the
steady state dynamics for the range $0\le\alpha\le70$.
There is no discernible difference between the bifurcation diagram of
the 8~mode nonlinear Galerkin approximation and the accurate
bifurcation diagram for this range of~$\alpha$.
Table~\ref{tablebifs2} lists the values of nonlinearity
parameter~$\alpha$ at which bifurcation points occur for the coarse
grid holistic models and the Galerkin approximations~\cite{Jolly90}.
The low mode first iterate nonlinear Galerkin approximations are
impressively accurate.

\begin{table}
\caption{$\alpha$ values at which bifurcation points occur for the
various coarse grid holistic models and low mode Galerkin
approximations}
\label{tablebifs2}
\begin{center}\footnotesize
\begin{tabular}{lcccccccc}
\hline
Approximation&$R_2b_1$&$R_2b_2$&$R_2b_3$&$R_2b_4$&$R_3t_1$&$R_3t_2$&$R_4b_1$&
$R_4q_1$\\ 
\hline
\multicolumn{4}{l}{Accurate $48$pts}\\ 
$6\mbox{th}$ order
&16.14&22.56&52.89&63.74&36.23&50.91&64.56&64.28\\ 
\hline
\multicolumn{4}{l}{Holistic~8 elements}
\\ 
$\Ord{\gamma^3,\alpha^2}$&14.64&20.36&39.34&44.96&$29.28^*$
&---&45.28&44.87\\ 
$\Ord{\gamma^3,\alpha^3}$&14.65&20.52&39.66&45.16&$29.33^*$
&---&45.47&44.96\\ 
$\Ord{\gamma^3,\alpha^4}$&14.65&20.53&39.72&45.21&$29.33^*$
&---&45.51&44.97\\ 
$\Ord{\gamma^4,\alpha^2}$&16.00&22.56&48.62&57.38&$34.73^*$&
---&57.89&57.49\\ 
$\Ord{\gamma^4,\alpha^3}$&16.00&22.56&48.25&56.84&$34.73^*$&
---&57.45&57.28\\ 
$\Ord{\gamma^4,\alpha^4}$&16.00&22.57&48.10&56.63&$34.73^*$&
---&57.30&57.21\\ 
$\Ord{\gamma^5,\alpha^2}$&16.13&22.72&51.54&61.54&$35.89^*$&
---&62.20&61.78\\ 
$\Ord{\gamma^5,\alpha^3}$&16.13&22.73&51.53&61.37&$35.91^*$&
---&62.04&61.70\\ 
$\Ord{\gamma^5,\alpha^4}$&16.13&22.73&51.60&61.38&$35.91^*$&
---&62.02&61.69\\ 
\hline
\multicolumn{4}{l}{Holistic~12 elements}\\ 
$\Ord{\gamma^3,\alpha^2}$&15.45&21.67&45.96&53.94&32.86
&45.98&54.49&54.17\\ 
$\Ord{\gamma^3,\alpha^3}$&15.45&21.69&46.05&54.00&32.87
&46.33&54.55&54.20\\ 
$\Ord{\gamma^4,\alpha^2}$&16.11&22.62&51.93&62.10&35.90
&50.92&62.83&62.52\\ 
$\Ord{\gamma^4,\alpha^3}$&16.11&22.62&51.94&62.10&35.90
&50.94&62.83&62.52\\ 
\hline
\multicolumn{4}{l}{Galerkin~\cite{Jolly90}} \\ 
3-m Euler--Galerkin&16.10&20.59&246.14&---&36.21&---&---&---\\ 
3-m Pseudo-stdy II&16.13&21.93&102.90&---&36.21&---&---&---\\ 
3-m Pseudo-stdy&16.13&22.01&93.91&---&36.24&63.91&---&---\\ 
12-m traditional&16.14&22.56&52.89&63.74&36.23&50.91&64.56&64.28\\ 
6-m traditional&16.14&22.55&52.72&63.28&36.23&46.85&64.00&64.00\\ 
3-m traditional&16.14&16.00&16.0??&16.0&36.00&36.0&---&---\\ 
\end{tabular}
\end{center}
\end{table}

This evidence suggests that the holistic models are competitive with
traditional Galerkin approximations, but that nonlinear Galerkin models
are significantly better.  However, recall that the holistic models are
based upon analysis of local dynamics and thus we expect them to be
more flexibly useful in applications than the global methods of these
Galerkin approximations.

\subsection{Coarse grids allow large time steps}
\label{S_ks_ss_ts}
A major benefit of accurate models on coarse grids is that larger time
steps are possible while maintaining numerical stability.
\S\ref{S_ks_ss_coarse} shows the remarkable accuracy of the
$\Ord{\gamma^5,\alpha^2}$ holistic model~(\ref{E_hol8dg3r1}) on a
coarse grid of 8~elements.
Here we investigate the maximum stable time step for \emph{explicit}
Runge--Kutta time integration on various holistic models---implicit
integration schemes are not considered.

\begin{figure}
\centering
\includegraphics[width=0.95\textwidth]{ks_stab_bif}
\caption{Bifurcation diagrams of (a)~$\Ord{\gamma^5,\alpha^2}$ holistic
model, with 8~elements on $[0,\pi]$, and (b)~2nd~order centered
difference approximation with 16~grid points on $[0,\pi]$.
Accurate bifurcation diagram is shown in grey.}
\label{ks_stab_bif}
\end{figure}

In particular we compare approximations of similar accuracy but
different grid resolutions to demonstrate the superior performance of
the holistic models.
For example, Figure~\ref{ks_stab_bif} compares the bifurcations
diagrams of the $\Ord{\gamma^5,\alpha^2}$ holistic model with
8~elements on $[0,\pi]$ and the 2nd~order centered difference
approximation~(\ref{E_con2_ks}) with 16~grid points on $[0,\pi]$.
The accurate bifurcation diagram is shown in grey.
See that the $\Ord{\gamma^5,\alpha^2}$ holistic model on the coarse grid
is of similar accuracy to the 2nd~order centered
difference approximation on the more refined grid.
Thus a reasonable comparison of computability is made using these two
schemes.

\begin{table}
\caption{Approximate maximum time steps for stability of 4th~order
Runge--Kutta scheme.}
\label{tablestabs}
\begin{center}\small
\begin{tabular}{llll}
\hline
Approximation&$\alpha=10$&$\alpha=20$&$\alpha=30$\\ 
\hline
\multicolumn{3}{l}{Holistic 8~elements}\\ 
$\Ord{\gamma^3,\alpha^2}$&.0011&.0014&.0017\\ 
$\Ord{\gamma^3,\alpha^3}$&.0011&.0014&.0017\\ 
$\Ord{\gamma^3,\alpha^4}$&.0011&.0014&.0017\\ 
$\Ord{\gamma^4,\alpha^2}$&.0006&.0007&.0008\\ 
$\Ord{\gamma^4,\alpha^3}$&.0006&.0007&.0008\\ 
$\Ord{\gamma^4,\alpha^4}$&.0006&.0007&.0008\\ 
$\Ord{\gamma^5,\alpha^2}$&.0005&.0005&.0006\\ 
$\Ord{\gamma^5,\alpha^3}$&.0005&.0005&.0006\\ 
$\Ord{\gamma^5,\alpha^4}$&.0005&.0005&.0006\\ 
\hline
\multicolumn{3}{l}{Centered 8~points}\\ 
2nd~order&.0011&.0012&---\\ 
4th~order&.0006&.0007&.0008\\ 
6th~order&.0005&.0005&.0006\\ 
\hline
\multicolumn{3}{l}{Centered 16~points}\\ 
2nd~order&.00006&.00006&.00006\\ 
\hline
\end{tabular}
\end{center}
\end{table}

Numerical experiments used the 4th~order Runge--Kutta scheme to
estimate the maximum stable time step for different holistic models and
centered difference approximations at various values of nonlinearity
parameter~$\alpha$.
Table~\ref{tablestabs} lists the approximate maximum time steps that
maintain numerical stability along both the negative unimodal branch at
$\alpha=10\,$, and the negative bimodal branch at $\alpha=20$ and
$\alpha=30\,$.
For the $\Ord{\gamma^5,\alpha^2}$ holistic model with 8~elements, the
maximum time step maintaining numerical stability is approximately
10~times larger than the corresponding time step for the 2nd~order
centered difference approximation with 16~grid points.
The $\Ord{\gamma^5,\alpha^2}$ holistic model requires approximately
3~times the number of floating point operations per grid value at each
time step compared to the 2nd~order centered difference approximation.
However, on a coarse grid of 16~points the 2nd~order centered
difference approximation must be applied at twice as many grid points.
Thus the $\Ord{\gamma^5,\alpha^2}$ holistic model can be integrated an
order of magnitude faster than the 2nd~order centered difference
approximation while maintaining similar accuracy.

Note: Table~\ref{tablestabs} shows that the higher order terms in the
nonlinearity~$\alpha$, generated by the holistic method, do not reduce
numerical stability.  Wider stencil holistic approximations reduce the
maximum stable time step somewhat, but so do the wider stencil
conventional centered difference approximations.  Thus, bear in mind
that we need to balance the accuracy gained by using higher order
approximation in~$\gamma$, that is, wider stencil approximations, with
the reduction in numerical stability and the increase in computation 
per grid value.

\section{Holistic models are accurate for time dependent phenomena}
\label{chap_ks_td}

The \KS\ equation~(\ref{EpdeKS}) has rich dynamics~\cite{Kuramoto78,
Jolly90, Kevrekidis90, Scovel88, Cross93, Holmes96, Dankowicz96,
Berkooz93}.
Having established the excellent performance of the holistic models in
reproducing the steady states of the \KS\ system in
Section~\ref{chap_ks_ss}, we now investigate the holistic models
performance at reproducing time dependent phenomena.
The \KS\ system exhibits complex time dependent behaviour such as limit
cycles, period doubling and spatio-temporal chaos.
This provides us with an example to explore the holistic approach to
modelling time dependent phenomena with relatively coarse
discretisations.

We restrict attention to $2\pi$~periodic solutions,
\begin{equation}
u(x,t)=u(x+2\pi,t)\,.
\label{E_odd_bc2}
\end{equation}
Initially we restrict further to solutions with odd symmetry, as in the
previous section, which exhibit, see Figure~\ref{bif6ord48}, Hopf
bifurcations to limit cycle solutions, and subsequent period doubling
bifurcations apparently leading to low-dimensional chaos~\cite{Jolly90,
Kevrekidis90, Scovel88}.
In \S\ref{S_numksper} we examine the dynamics of the holistic models on
coarse grids through the eigenvalues of the models near the steady
states.
For example, we see that the $\Ord{\gamma^5,\alpha^2}$ holistic model
reproduces much of the eigenvalue information for $0\le\alpha\le70$ on
a coarse grid of 8~elements.
In~\S\ref{S_ks_td_hb} we explore the bifurcation diagrams near the
first Hopf bifurcation and capture the stable limit cycles and period
doubling sequence.
The holistic models more accurately model the dynamics than centered
difference approximations of equal stencil width.
Subsequently we just require spatial periodicity whence stable
travelling wave appear followed by, at higher values of nonlinearity
parameter~$\alpha$, more complex spatio-temporal chaos as investigated
by Holmes, Lumley \& Berkooz~\cite{Holmes96} and
Dankowicz et~al.~\cite{Dankowicz96}.
In~\S\ref{S_ks_td_no} we find the holistic discretisations more
accurately model the amplitude and wave speed of travelling wave
solutions, and predict better space time plots and time averaged power
spectra, than corresponding the centered difference approximations.

\subsection{Dynamics near the steady states are reproduced}
\label{S_numksper}
Consider the eigenvalues of the \KS{} system~(\ref{EpdeKS}) linearised
about the the steady states and restricted to odd symmetry.
Accurate modelling of the eigenvalues near the steady states is a
necessary condition for the accurate modelling of the dynamics.
We look at two views of the eigenvalues: first, their value on the
negative bimodal branch; and second a more qualitative plot of their
values on the entire bifurcation diagram for nonlinearity parameter
$0\le\alpha\le70$\,.

\begin{figure}
\centering
\includegraphics[width=\textwidth]{eigshol2}
\caption{The four largest (least negative) eigenvalues along the stable
bimodal branch for the (a)~$\Ord{\gamma^3,\alpha^2}$,
(b)~$\Ord{\gamma^4,\alpha^2}$, (c)~$\Ord{\gamma^5,\alpha^2}$ holistic
models shown in green for 8~elements on $[0,\pi]$.  The accurate
eigenvalues are shown in blue.}
\label{eigshol}
\end{figure}

\paragraph{Compare eigenvalues along the bimodal branch}
We investigate dynamics near the stable negative bimodal branch.
Consider the real part of the four largest (least negative real part)
eigenvalues for low order holistic models and compare to explicit
centered difference approximations on a coarse grid of 8~elements
on~$[0,\pi]$.
Figure~\ref{eigshol} shows the four largest eigenvalues for the
$\Ord{\gamma^3,\alpha^2}$~(\ref{E_hol8dg1r1}),
$\Ord{\gamma^4,\alpha^2}$~(\ref{E_hol8dg2r1}) and
$\Ord{\gamma^5,\alpha^2}$~(\ref{E_hol8dg3r1}) holistic models in green
and the accurate solution in blue.\footnote{As in
Section~\ref{chap_ks_ss}, the accurate reference for solutions is found
using a 6th~order centered difference approximation with 48~grid points
on~$[0,\pi]$.} Recall the $\Ord{\gamma^3,\alpha^2}$,
$\Ord{\gamma^4,\alpha^2}$ and $\Ord{\gamma^5,\alpha^2}$ holistic models
have 5~point, 7~point and 9~point stencils, respectively.
Figure~\ref{eigshol}c, shows the four largest eigenvalues for the
$\Ord{\gamma^5,\alpha^2}$ holistic model closely matches the accurate
solution over this range of nonlinearity parameter~$\alpha$.

\begin{figure}
\centering
\includegraphics[width=\textwidth]{eigscon2}
\caption{The four largest eigenvalues along the stable bimodal branch
for the (a)~2nd~order, (b)~4th~order, (c)~6th~order centered difference
approximations shown in magenta for 8~grid points on~$[0,\pi]$.
The accurate eigenvalues are shown in blue.}
\label{eigscon}
\end{figure}

Similarly, Figure~\ref{eigscon} shows the four largest eigenvalues for
the 2nd~order~(\ref{E_con2_ks}), 4th~order~(\ref{E_con4_ks}) and
6th~order~(\ref{E_con6_ks}) centered difference approximations in
magenta on the same coarse grid.
The centered difference approximations shown here are of equal stencil
width to the corresponding holistic models in Figure~\ref{eigshol}.
Figure~\ref{eigscon}a, shows the 2nd~order centered difference barely
approximates the behaviour of the stable bimodal branch for
$\alpha<20$\,.
Even the 6th~order centered difference
approximation,~Figure~\ref{eigscon}c, is inferior to the
$\Ord{\gamma^4,\alpha^2}$ holistic model for $\alpha>30$\,.
This is despite the 6th~order centered difference model having a wider
stencil of 9~points compared to the 7~point stencil of the
$\Ord{\gamma^4,\alpha^2}$ holistic model.
Figures \ref{eigshol}~and~\ref{eigscon} show the low order holistic
models are superior to the corresponding centered difference
approximations for reproducing the dynamics near the stable bimodal
branch.

\begin{figure}
\centering
\includegraphics[width=\textwidth]{hol8eigs2}
\caption{Bifurcation diagram of the $\Ord{\gamma^5,\alpha^2}$ holistic
model with 8~elements and odd symmetry on~$[0,\pi]$, depicting the real parts of the
8~largest (least negative) eigenvalues colour coded according to the
colour bar shown.  } \label{hol8eigs}
\end{figure}

\begin{figure}
\centering
\includegraphics[width=0.95\textwidth]{con6ord48eigs2}
\caption{Bifurcation diagram of the accurate \KS\ system, depicting the
real parts of the 8~largest (least negative) eigenvalues, colour coded
according to the colour bar shown.} \label{con6ord48eigs}
\end{figure}

\paragraph{Compare eigenvalues across the bifurcation diagram}
Here we explore a new view of the earlier bifurcation diagrams that
additionally depicts the \emph{real part} of the 8~largest (least
negative) eigenvalues by colour.
Compare the eigenvalues of the
$\Ord{\gamma^5,\alpha^2}$~(\ref{E_hol8dg3r1}) holistic model, see
Figure~\ref{hol8eigs}, on the coarse grid of 8~elements on~$[0,\pi]$ to
accurate ones for the \KS\ system, see Figure~\ref{con6ord48eigs}, over
the nonlinearity parameter $0\le\alpha\le70$\,.
The magnitude of the real part of the eigenvalues is colour coded
according to the colour bar shown on the right of the bifurcation
diagram; the least negative eigenvalues are plotted above the more
negative to give a small band of colour for each branch of steady
states at each parameter value.
Similarly to the bifurcation diagrams shown in
Section~\ref{chap_ks_ss}, the open squares denote bifurcation points
and the black squares denote Hopf bifurcations.
Figure~\ref{hol8eigs}, when compared to Figure~\ref{con6ord48eigs},
shows that in addition to reproducing the stability of the accurate
\KS\ system for $0\le\alpha\le70$ as discussed
in~\S\ref{S_ks_ss_coarse}, the $\Ord{\gamma^5,\alpha^2}$ holistic model
reproduces well the eigenvalues for most of this range of nonlinearity
parameter~$\alpha$.
This accurate modelling of the eigenvalues is evidence of accurate
modelling of the\KS\ dynamics, at least near the steady states.

\subsection{Extend the Hopf bifurcations}
\label{S_ks_td_hb}

Hopf bifurcations give rise to time periodic solutions  (limit cycles).
We explore the predictions of the various models to see how well they
capture these strongly time dependent phenomena.

Here we investigate the bifurcation diagrams obtained by extending the
first Hopf bifurcation, at $\alpha=30.345$\,, on the positive bimodal
branch and the period doubling sequence that ensues.
We compare the bifurcation diagrams of low order holistic models to
explicit centered difference models on a coarse grid of 8~elements on
$[0,\pi]$ to the accurate bifurcation diagram of the \KS\ system.
% We also investigate the period doubling sequence which continues from
% the first Hopf bifurcation and compare the
% $\Ord{\gamma^5,\alpha^2}$~(\ref{E_hol8dg3r1}) holistic models
% approximation of the period doubling sequence to the 6th~order centered
% difference approximation~(\ref{E_con6_ks}) and the accurate period
% doubling sequence.
Trajectories in the period doubling sequence are reported and compared
by MacKenzie~\cite{MacKenzie05}.
As before, the holistic models outperform the corresponding centered
difference approximations.

\begin{figure}
\centering
\includegraphics[width=0.95\textwidth]{per_bifhol_8}
\caption{Bifurcation diagrams near the first Hopf bifurcation for
(a)~$\Ord{\gamma^3,\alpha^2}$, (b)~$\Ord{\gamma^4,\alpha^2}$,
(c)~$\Ord{\gamma^5,\alpha^2}$ holistic models with 8~elements on
$[0,\pi]$ and (d)~an accurate bifurcation diagram.  Stable limit cycles
are shown in light blue and unstable limit cycles are shown in orange.}
\label{per_bifhol_8}
\end{figure}

\paragraph{Investigate the first Hopf bifurcation}
We now investigate the holistic models near the first Hopf bifurcation
on the positive bimodal branch, labeled~$\text{\textsc{hb}}_1$, with a coarse grid of
8~elements on $[0,\pi]$.
Figure~\ref{per_bifhol_8} shows the bifurcation diagrams of the low
order holistic models and the accurate bifurcation diagram near the
first Hopf bifurcation.
The stable limit cycles (light blue) that continue from this
bifurcation point undertake a period doubling sequence commencing at a
point labeled~$\text{\textsc{pd}}$ (yellow square).
The pair of unstable limit cycles born at~$\text{\textsc{pd}}$ give rise to the period
doubling sequence leading to chaos.

The accurate bifurcation diagram shown in Figure~\ref{per_bifhol_8}d,
is produced using a 6th~order centered difference approximation with
24~grid points on~$[0,\pi]$.
The accurate bifurcation diagram shown is identical to the bifurcation
diagram for the same range of~$\alpha$ produced by Jolly
{et~al.}~\cite{Jolly90}.
Figure~\ref{per_bifhol_8}a, shows that even the lowest order
$\Ord{\gamma^3,\alpha^2}$~(\ref{E_hol8dg1r1}) holistic model reproduces
the the first Hopf bifurcation and finds the period doubling point on
this coarse grid of 8~elements\footnote{Figure~\ref{per_bifhol_8}a
displays the bifurcation diagram for $25\le\alpha\le32$ compared to
$30\le\alpha\le37$ for the other diagrams.
Since the first Hopf bifurcation for the $\Ord{\gamma^3,\alpha^2}$
holistic model occurs at $\alpha=25.595$ the bifurcation diagram is
shifted to contain the important dynamics.
} In comparison, the corresponding 2nd~order centered difference
approximation does not even have the first Hopf bifurcation, see
Figure~\ref{con_ks_8bif_all}b.

\begin{figure}
\centering
\includegraphics[width=0.95\textwidth]{per_bifcon_8}
\caption{Bifurcation diagrams near the first Hopf bifurcation for (a)
4th~order, (b) 6th~order centered difference approximations with 8~grid
points on $[0,\pi]$.
Stable limit cycles are shown in light blue and unstable limit cycles
are shown in orange.}
\label{per_bifcon_8}
\end{figure}
Figure~\ref{per_bifhol_8}b,c, show that higher order holistic models
accurately model the first Hopf bifurcation and the resulting stable
and unstable limit cycles.
The accuracy of the $\Ord{\gamma^5,\alpha^2}$~(\ref{E_hol8dg3r1})
holistic model for reproducing these periodic solutions of the \KS{}
system is remarkable on this coarse grid.
Figure~\ref{per_bifcon_8} shows the corresponding bifurcation diagrams
for the 4th~order and 6th~order centered difference approximations with
8~grid points on $[0,\pi]$.
Compare Figure~\ref{per_bifcon_8} and Figure~\ref{per_bifhol_8} to see
that the 6th~order centered difference approximation which has a nine
point stencil does not perform as well as the
$\Ord{\gamma^4,\alpha^2}$~(\ref{E_hol8dg2r1}) holistic model which has
a 7~point stencil.
Figure~\ref{per_bifhol_8}b,c, show that higher order holistic models
more accurately model the first Hopf bifurcation and the resulting
stable and unstable limit cycles.
The accuracy of the $\Ord{\gamma^5,\alpha^2}$~(\ref{E_hol8dg3r1})
holistic model for reproducing these periodic solutions of the \KS{}
system is remarkable on this coarse grid.
Table~\ref{tableper} shows the parameter values $\alpha$ for the Hopf
bifurcations, $\text{\textsc{hb}}_1$ and the initial period doubling point~$\text{\textsc{pd}}$.
See that both the $\Ord{\gamma^4,\alpha^2}$ and
$\Ord{\gamma^5,\alpha^2}$ holistic models are more accurate than the
4th~order and 6th~order centered difference approximations in
reproducing the first Hopf bifurcation and the resulting period
doubling point.

\begin{table}
\caption{Nonlinearity parameter~$\alpha$ values for the first Hopf
bifurcation point~$\text{\textsc{hb}}_1$ and resulting period doubling point~$\text{\textsc{pd}}$.}
\label{tableper}
\begin{center}\small
\begin{tabular}{lcc}
\hline
Approximation&$\text{\textsc{hb}}_1$&$\text{\textsc{pd}}$\\ 
\hline
Holistic 8~elements\\ 
$\Ord{\gamma^3,\alpha^2}$&25.60&27.22\\ 
$\Ord{\gamma^4,\alpha^2}$&30.04&32.03\\ 
$\Ord{\gamma^5,\alpha^2}$&30.66&32.95\\ 
\hline
Centered 8~points\\ 
2nd~order&---&---\\ 
$4$th order&27.91 &29.57 \\ 
$6$th order&29.11 &31.40 \\ 
\hline
Accurate&30.35&32.97\\ 
\hline
\end{tabular}
\end{center}
\end{table}

\subsection{Dynamics of periodic patterns without odd symmetry}
\label{S_ks_td_no}

Consider the \KS\ system~(\ref{EpdeKS}) with solutions that are
spatially periodic~(\ref{E_odd_bc2}).
We remove the requirement for odd symmetry.
Consequently, we now explore travelling wave like solutions at low
nonlinearity~$\alpha$.
Also, we investigate the spatio-temporal chaos that occurs at
higher~$\alpha$.

\begin{figure}
\centering
\includegraphics[width=0.95\textwidth]{kuraf_t_a5_8}
\caption{${\alpha=5}$\,: wave-like solutions  at
$t=0,0.2,0.4,0.6,0.8,1\,$ for the $\Ord{\gamma^3,\alpha^2}$ holistic
model shown in green and the 2nd~order centered difference
approximation in magenta on a coarse grids of 8~elements on $[0,2\pi]$.
The accurate solution is shown in blue.}
\label{kuraf_t_a5_8}
\end{figure}

\paragraph{Good performance for holistic models at low $\alpha$}
\label{S_ks_td_low_a}
Consider the holistic models of the \KS{} system~(\ref{EpdeKS}) and
(\ref{E_odd_bc2}) for nonlinearity parameter $\alpha=5$ and $\alpha=10$
on coarse grids of 8~elements on~$[0,2\pi]$.
Figure~\ref{kuraf_t_a5_8} shows solutions obtained from the lowest
order $\Ord{\gamma^3,\alpha^2}$~(\ref{E_hol8dg1r1}) holistic model for
$\alpha=5$ in green, the accurate solution in blue and the
corresponding 2nd~order centered difference
approximation~(\ref{E_con2_ks}) with 8~points on~$[0,2\pi]$, in magenta.
The solutions are shown at time slices $t=0,0.2,0.4,0.6,0.8,1\,$,
starting from the half-wave initial condition of $u(x,0)=|\sin(x/2)|$\,.
See the $\Ord{\gamma^3,\alpha^2}$ holistic model is superior to the
2nd~order centered difference approximation on this coarse grid.
In particular, the amplitude of the evolving wave-like solution and the
wave speed are more accurately reproduced by the
$\Ord{\gamma^3,\alpha^2}$ holistic model for $\alpha=5$\,.

\begin{figure}
\centering
\includegraphics[width=0.95\textwidth]{kuraf_t_a10_8_2}
\caption{${\alpha=10}$\,: wave-like solutions at
$t=0,0.2,0.4,0.6,0.8,1\,$ for the $\Ord{\gamma^3,\alpha^2}$,
$\Ord{\gamma^4,\alpha^2}$ and $\Ord{\gamma^5,\alpha^2}$ holistic models
shown in green, light green and light blue respectively and the
6th~order centered difference approximation shown in red on a coarse
grids of 8~elements on $[0,2\pi]$.  The accurate solution is shown in
blue.}
\label{kuraf_t_a10_8}
\end{figure}

Similarly, Figure~\ref{kuraf_t_a10_8} shows the same time slices for
larger $\alpha=10$\,.
The $\Ord{\gamma^3,\alpha^2}$~(\ref{E_hol8dg1r1}) holistic model is
shown in green, the $\Ord{\gamma^4,\alpha^2}$~(\ref{E_hol8dg2r1}) model
is shown in light green and the
$\Ord{\gamma^5,\alpha^2}$~(\ref{E_hol8dg3r1}) holistic model in light
blue for this coarse grid of 8~elements on~$[0,2\pi]$.
For this~$\alpha$ the 2nd~order~(\ref{E_con2_ks}) and
4th~order~(\ref{E_con4_ks}) centered difference approximations do not
generate a wave-like solution at all.
However, the 6th~order centered difference
approximation~(\ref{E_con6_ks}) does produce the travelling wave-like
solution shown in red.
The $\Ord{\gamma^3,\alpha^2}$ holistic model (green) is the least
accurate on this coarse grid but it does reproduce a stable solution on
this coarse grid for only a 5~point stencil approximation.
The $\Ord{\gamma^4,\alpha^2}$ holistic model (light green) more
accurately models the amplitude of the solution compared to the
6th~order centered difference approximation despite having a smaller
stencil width.
The $\Ord{\gamma^5,\alpha^2}$ holistic model is the most accurate at
reproducing both the amplitude and the wave speed of the stable
wave-like solution for $\alpha=10$ on this coarse grid of 8~elements.

\paragraph{Good performance for more complex behaviour}
For higher values of nonlinearity parameter~$\alpha$ for which the \KS\
system exhibits more complex behaviour, including spatio-temporal
chaos, it is more useful to compare time averaged power spectra
rather than particular travelling waves.
Here we investigate the performance of the holistic models on coarse
grids for $\alpha=20$~and~$50$ using the example of the
$\Ord{\gamma^5,\alpha^2}$~(\ref{E_hol8dg3r1}) holistic model on
relatively coarse grids, and we compare it with the 6th~order centered
difference approximation which is of equal stencil width.
The focus here is to show the improved performance of the holistic
models for ranges of parameter~$\alpha$ that contain more complex time
dependent behaviour.
We also expect a corresponding improvement for the other holistic
models but this is not investigated here.
Further, we also compare the $\Ord{\gamma^5,\alpha^2}$ holistic model
on coarse grids to the 2nd~order centered difference approximations of
similar accuracy.
We find the $\Ord{\gamma^5,\alpha^2}$ holistic model, but with
approximately 1/3 of the grid points, has comparable accuracy to
2nd~order centered difference approximations.
 
\begin{figure}
\centering
\includegraphics[width=\textwidth]{ks_td_st_a20_1}
\caption{{$\alpha=20$}\,: space time plots for (a)~the
$\Ord{\gamma^5,\alpha^2}$ holistic model with 12~elements on
$[0,2\pi]$, (b)~6th~order centered difference approximation with
12~grid points on $[0,2\pi]$ and (c)~the accurate solution}
\label{ks_td_st_a20_1}
\end{figure}

Figure~\ref{ks_td_st_a20_1} shows space time plots of (a)~the
$\Ord{\gamma^5,\alpha^2}$~holistic model with 12~elements on
$[0,2\pi]$, (b)~the 6th~order centered difference approximation with
12~grid points on $[0,2\pi]$ and (c)~the accurate
solution.\footnote{The accurate solutions plotted in this section are
computed using a 6th~order centered difference approximation and 256
grid points on the interval $[0,2\pi]$.
This is sufficient grid resolution to capture the important dynamics of
the \KS\ system for the values of~$\alpha$ investigated here.} See the
$\Ord{\gamma^5,\alpha^2}$~holistic model reproduces much of the complex
structure of the accurate solution for $\alpha=20$ with 12~elements.
Figure~\ref{ks_td_st_a20_1}b, shows the 6th~order centered difference
approximation incorrectly finds a periodic solution after approximately
$t=0.2$\,.

\begin{figure}
\centering
\includegraphics[width=\textwidth]{ks_td_ps_a20_1}
\caption{{$\alpha=20$}\,: time averaged power spectra  for the
$\Ord{\gamma^5,\alpha^2}$ holistic model with 12~elements on $[0,2\pi]$
shown in light blue, and the 6th~order centered difference
approximation in red for (a)~12~grid points on $[0,2\pi]$ and
(b)~16~grid points on $[0,2\pi]$.  The accurate power spectrum is shown
in blue.}
\label{ks_td_ps_a20_1}
\end{figure}

Since the \KS\ system at nonlinearity parameter $\alpha=20$ exhibits
more complex time dependent behaviour than simple limit cycles, we
compare time averaged power spectra, denoted here by~$S(k)$ for
wavenumber~$k$.
Figure~\ref{ks_td_ps_a20_1}a, shows a log-log plot of the time average
power spectra of the $\Ord{\gamma^5,\alpha^2}$ holistic model in
light blue and the 6th~order centered difference approximation on a
coarse grid of 12~elements on $[0,2\pi]$ in red.
The accurate power spectrum is shown in blue.
For this coarse grid of only 12~elements only 5~wavenumbers are
relevant, as displayed.
See that the $\Ord{\gamma^5,\alpha^2}$ holistic model is superior to
the 6th~order centered difference approximation on this coarse grid of
12~elements.
Figure~\ref{ks_td_ps_a20_1}b, compares the time average power spectrum
of the $\Ord{\gamma^5,\alpha^2}$ holistic model in light blue with
12~elements and the 6th~order centered difference approximation with
16~grid points.
The $\Ord{\gamma^5,\alpha^2}$ holistic model achieves similar accuracy
on a coarser grid.
 
\begin{figure}
\centering
\includegraphics[width=\textwidth]{ks_td_ps_a20_3}
\caption{{$\alpha=20$}\,: time averaged power spectra for  for the
$\Ord{\gamma^5,\alpha^2}$ holistic model with 12~elements on $[0,2\pi]$
shown in light blue, and the 2nd~order centered difference
approximation in magenta for (a)~24~grid points on $[0,2\pi]$ and (b)~36~grid points on $[0,2\pi]$.
The accurate spectrum is shown in blue.}
\label{ks_td_ps_a20_3}
\end{figure}

The power spectra of the $\Ord{\gamma^5,\alpha^2}$ holistic model on a
coarse grid of 12~elements and the 2nd~order centered difference
approximation on the more refined grids of 24~and~36 points are shown
in Figures~\ref{ks_td_ps_a20_1}a,b respectively.
See a refined grid of 36~points is needed achieve similar accuracy to
the $\Ord{\gamma^5,\alpha^2}$ holistic model on a coarse grid of
12~elements on~$[0,2\pi]$.
That is, through its subgrid scale modeling, the holistic model
achieves similar accuracy with one-third the dimensionality.

\begin{figure}
\centering
\includegraphics[width=\textwidth]{ks_td_st_a50_1}
\caption{{$\alpha=50$}\,: space time plots for (a)~the
$\Ord{\gamma^5,\alpha^2}$ holistic model with 24~elements
on~$[0,2\pi]$, (b)~6th~order centered difference approximation with
24~grid points on~$[0,2\pi]$ and (c)~the accurate solution.}
\label{ks_td_st_a50_1}
\end{figure}

For nonlinearity parameter $\alpha=50$ the \KS\ system exhibits even
more complex behaviour, see the space time plots in
Figure~\ref{ks_td_st_a50_1}.
The $\Ord{\gamma^5,\alpha^2}$ holistic model more accurately reproduces
the \KS\ system than the 6th~order centered difference approximation on
this coarse grid of 24~elements.
On this coarse grid the 6th~order centered difference approximation
shown in Figure~\ref{ks_td_st_a50_1}b, exhibits a periodic solution
after time $t\approx 0.1$ which does not match the irregular behaviour
seen in the accurate solution and the $\Ord{\gamma^5,\alpha^2}$
holistic model.

\begin{figure}
\centering
\includegraphics[width=\textwidth]{ks_td_ps_a50_1}
\caption{{$\alpha=50$}\,: time averaged power spectra for the
$\Ord{\gamma^5,\alpha^2}$ holistic model with 24~elements on~$[0,2\pi]$
shown in light blue, and the 6th~order centered difference
approximation in red for (a)~24~grid points on~$[0,2\pi]$ and
(b)~32~grid points on~$[0,2\pi]$.  The accurate spectrum is shown in
blue.}
\label{ks_td_ps_a50_1}
\end{figure}

We again examine time averaged power spectra to further investigate the
performance of the $\Ord{\gamma^5,\alpha^2}$ holistic model at this
relatively large parameter value of $\alpha=50$.
Figure~\ref{ks_td_ps_a50_1} compares the time averaged power spectrum
of the $\Ord{\gamma^5,\alpha^2}$ holistic model in light blue on a
coarse grid of 24~elements on $[0,2\pi]$ to the 6th~order centered
difference approximation in red, for (a)~24~grid points and (b)~32~grid
points on $[0,2\pi]$.
The 6th~order centered difference approximation with 32~grid points has
similar accuracy to the $\Ord{\gamma^5,\alpha^2}$ holistic model on a
coarse grid of just 24~elements for $\alpha=50$\,.

This investigation of the $\Ord{\gamma^5,\alpha^2}$ holistic model on
coarse grids for $\alpha=20$ and~$50$ shows it reproduces similar
accuracy to the 2nd~order centered difference approximation on a coarse
grid of approximately~1/3 the resolution, and similar accuracy to the
6th~order centered difference approximation on grids of approximately
3/4 the resolution.  MacKenzie~\cite{MacKenzie05} reports that even
at $\alpha=200$ the holistic model  qualitatively well captures the
dynamics of the \KS\ \pde.
This increased accuracy on coarse grids allows larger time steps for
explicit time integration schemes, as discussed in \S\ref{S_ks_ss_ts}.

\section{Conclusion}

Holistic discretisation~\cite{Roberts98a} is straightforwardly extended
to fourth order dissipative \pde s through the example of the \KS\
equation~\cite{MacKenzie00a}.
We divide the domain into elements by introducing artificial internal
boundary conditions~(\S\ref{S_KSho}) which isolate the elements when
$\gamma=0$ but when $\gamma=1$ they fully couple the elements to
recover the \KS\ equation.
Then centre manifold theory supports the
discretisation,  see \S\ref{S_KScm}.
The holistic models listed in~\S\ref{S_ks_rel} have a dual
justification~(\S\ref{S_ks_epde}): not only are they supported by
centre manifold theory for finite element size~$h$, the \ibc{} are
specially crafted~\cite{Roberts00a} so the models are also consistent
with the \KS{} equation as the grid spacing $h\rightarrow0\,$.

No formal error bounds currently exist for the holistic method; the
difficulty is that the models are based at $\gamma=0$ but are evaluated
at finite $\gamma=1$\,.  Instead we present a detailed numerical
investigation of the holistic models of the steady
states~(Section~\ref{chap_ks_ss}) and time dependent
solutions~(Section~\ref{chap_ks_td}) of the \KS\ on coarse grids.

We compared, in \S\ref{S_ks_ss_ts}, the accuracy of different
approximations in predicting steady states on different grid
resolutions.
The holistic $\Ord{\gamma^5,\alpha^2}$ approximation on a grid of
8~elements has similar accuracy to a 2nd~order centered difference
approximation on a grid of 16~points.
Consequently the holistic model allows a maximum time step which is an
order of magnitude longer than that of the explicit centered difference
approximation of similar accuracy, while maintaining numerical
stability.
The accuracy of the holistic approximations to the \KS{} equation on
coarse grids and subsequent improved performance justifies further
application of the holistic method and future investigation of the
approach.

The holistic models on coarse grids also modelled well time dependent
phenomena of the \KS\ system.
In particular, in \S\ref{S_numksper} we saw the holistic models more
accurately model the eigenvalues near the steady states of the first
form of the \KS\ system compared to explicit centered difference
approximations of equal stencil widths.
The coarse grid holistic models also more accurately model
the first Hopf bifurcation and the resulting period doubling sequence,
see \S\ref{S_ks_td_hb}.
Further, in comparison with explicit centered difference models, in
\S\ref{S_ks_td_no}, we saw good performance for higher values of the
nonlinearity parameter~$\alpha$ and more accurate predictions of time
averaged power spectra: the $\Ord{\gamma^5,\alpha^2}$ holistic model
achieves similar accuracy to the 2nd~order and 6th~order centered
difference approximations on approximately 1/3 and 3/4 of the grid
resolutions respectively.

This good performance of the holistic models for accurately reproducing
both the steady states and the time dependent phenomena of the \KS\
system is good evidence that the holistic approach is a powerful method
for discretising dissipative \pde{}s on coarse grids.

\appendix
\section{Computer algebra derives the discretisation}
\label{app_a}

{\small
\verbatimlisting{pde_dbc.red}
}

\bibliographystyle{plain}
\bibliography{ajr,new,bib}

\end{document}